\documentclass[a4paper,11pt]{amsart}
\usepackage{amssymb,amscd,amsfonts,amsbsy,nicefrac}
\usepackage{amsmath,amsthm,amssymb,amsfonts,enumitem,verbatim}
\usepackage{ bbold }
\usepackage[hidelinks]{hyperref}
\usepackage{amsrefs}
\usepackage{xcolor}
\usepackage{mathrsfs}
\usepackage{graphicx}

\newcommand{\nocomment}[1]{}

\calclayout

\usepackage{enumitem}
\usepackage{parskip}
\usepackage{unicodetex}

\makeatletter
\g@addto@macro\bfseries{\boldmath}
\makeatother

\makeatletter
\def\subsubsection{\@startsection{subsubsection}{3}%
\z@{.5\linespacing\@plus.7\linespacing}{-.5em}%
{\normalfont\bfseries}}
\makeatother

\newtheorem{thm}{Theorem}[section]
\newtheorem{lem}[thm]{Lemma}
\newtheorem{lemma}[thm]{Lemma}

\newtheorem{prop}[thm]{Proposition}
\newtheorem{defn}[thm]{Definition}
\newtheorem{rem}[thm]{Remark}

\def\vth{\vartheta}
\def\th{\theta}
\def\psiF{\psi}
\def\phiF{\phi}
\def\omegaF{\omega}

\def\Im{{\rm Im}\,}
\def\R{\mathbb{R}}

\def\eps{\varepsilon}

\def\2q{{\frac{2}{|B|}}}

\newcommand{\N}{\mathbb{N}}

\newcommand{\bmo}{\mathrm{bmo}}
\newcommand{\BMO}{\mathrm{BMO}}

\newcommand{\Z}{\mathbb{Z}}

\newcommand{\jap}[1]{\langle{#1}\rangle}

\newcommand{\supp}{\mbox{supp}\,}

\renewcommand{\leq}{\leqslant}
\renewcommand{\geq}{\geqslant}

\newcommand{\phase}{\varphi}
\newcommand{\Phase}{\Phi}

\newcommand{\Mod}[1]{\ (\mathrm{mod}\ #1)}

\newcommand{\m}[1]{\begin{equation*}
#1
\end{equation*}}
\newcommand{\nm}[2]{\begin{equation}\label{#1}
#2
\end{equation}}
\newcommand{\ma}[1]{\begin{align*}
#1
\end{align*}}
\newcommand{\nma}[2]{\begin{equation}\label{#1}
\begin{aligned}
#2
\end{aligned}
\end{equation}}

\newcommand{\abs}[1]{\left|#1\right|}
\newcommand{\set}[1]{\left\{#1\right\}}
\newcommand{\brkt}[1]{\left(#1\right)}
\newcommand{\dd}{\,\mathrm{d}}
\newcommand{\dddd}{\mathrm{d}}
\newcommand{\ddd}{\,\text{\rm{\mbox{\dj}}}}

\newcommand{\norm}[1]{\left\Vert#1\right\Vert}

\begin{document}

\title[Estimates and multilinear oscillatory integrals] {Multilinear oscillatory integrals and estimates for coupled systems of dispersive PDEs}

\author[A.~Bergfeldt]{Aksel Bergfeldt}
\address{A.~Bergfeldt, Department of Mathematics, Uppsala University, SE-751 06 Uppsala, Sweden}
\email{aksel.bergfeldt@math.uu.se}
\author[S.~Rodr\'iguez-L\'opez]{Salvador Rodr\'iguez-L\'opez}
\address{S.~Rodr\'iguez-L\'opez, Department of Mathematics, Stockholm University, SE-106 91 Stockholm, Sweden}
\email{s.rodriguez-lopez@math.su.se}
\author[D.~Rule]{David Rule}
\address{D.~Rule, Department of Mathematics, Link\"oping University, SE-581 83 Link\"oping, Sweden}
\email{david.rule@liu.se}
\author[W.~Staubach]{Wolfgang Staubach}
\address{W.~Staubach, Department of Mathematics, Uppsala University, SE-751 06 Uppsala, Sweden}
\email{wulf@math.uu.se}
\thanks{The second author has been partially supported by the Grant PID2020-113048GB-I00. The third author was partially supported by the Research School in Interdisciplinary Mathematics at Link\"oping University. }

\begin{abstract}
We establish sharp global regularity of a class of multilinear oscillatory integral operators that are associated to nonlinear dispersive equations with both Banach and quasi-Banach target spaces. As a consequence we also prove the (local in time) continuous dependence on the initial data for solutions of a large class of coupled systems of dispersive partial differential equations.  
\end{abstract}

\maketitle

\section{Introduction}\label{sec:intro}

In this paper, we consider the regularity of multilinear oscillatory integral operators (multilinear OIOs for short) that are associated to nonlinear dispersive equations in the realm of Banach and quasi-Banach function spaces. Examples include non-linear water-wave and capillary wave equations, nonlinear wave and Klein--Gordon equations, the nonlinear Schr\"odinger equations, the Korteweg--de\,Vries-type equations, and higher order nonlinear dispersive equations. To achieve this, we develop a fairly general and complete framework  for the investigation of the regularity of a class of multilinear oscillatory integral operators with smooth amplitudes.\\

The literature on multilinear oscillatory integrals is by now quite vast. However if we confine ourselves to those operators that appear in connection to non-linear PDEs, then one only has a handful of optimal (i.\,e.\@{} endpoint) results. These are:

\begin{enumerate}[label=\textit{\roman*})]
    \item Sharp global regularity of bilinear and multilinear oscillatory integral operators (where the phase function of the operator is homogeneous of degree one) with Banach target spaces see \cite{Monster} and \cite{RRS2}.
    
    \item Sharp global regularity of bilinear oscillatory integral operators that also include operators with quadratic behaviour in their phase functions, and with Lebesgue-type targets in the Banach scales, see the work of F. Bernicot and P. Germain \cite{BernicotGermain1}. Also some sharp global results for certain multilinear operators with quadratic phase functions (different from those considered in \cite{BernicotGermain1}) were obtained in \cite{BS}. 
\end{enumerate}
The main contributions of this paper can be briefly summarised as follows:
\begin{itemize}
    \item From the point of view of nonlinear PDEs, we prove regularity results that can be used in understanding the interaction of free waves in coupled systems of PDE's which can in turn be used to  understand more complicated  nonlinear problems.
    \item From the point of view of Fourier analysis, we extend the regularity of multilinear OIOs with homogeneous of degree one phase functions to the case of operators with inhomogeneous phase functions, and with target spaces that include quasi-Banach as well as Banach Hardy spaces. In this context and for the given multiplier operators at hand, our results are optimal.
\end{itemize}

In proving boundedness results with quasi-Banach target spaces, although an appropriate frequency-space decomposition of the operators into various frequency regimes is available (see e.\,g.\@ \cite{RRS2}), the classical duality method of R. Coifman and Y. Meyer \cite{CM}, is no longer applicable. Therefore, in this paper we introduce an approach based on

\begin{enumerate}
 \item vector-valued inequalities,
    \item maximal functions of Hardy--Littlewood, Peetre and Park,
    \item estimates for linear oscillatory integral operators,
\end{enumerate}
which enable us to prove the desired boundedness results.

\subsection{Some results concerning boundedness of multilinear oscillatory integral operators}
We start by giving an overview of the previously known regularity results for OIOs, which are relevant to operators that are considered here.

\begin{defn}\label{defn:amplitudes}
For integers $n,N\geq 1$ and $m \in \R$, the set of \emph{(multilinear) amplitudes}  $S^m(n,N)$ is the set of functions $\sigma \in \mathcal{C}^\infty (\R^n \times \R^{nN})$ that satisfy
\m{
|\partial_{\Xi}^{\alpha}\partial_{x}^{\beta} \sigma(x,\Xi)|\leq C_{\alpha,\beta} \langle \Xi\rangle^{m-|\alpha|},
}
for all multi-indices $\alpha$ and $\beta$. Here and in what follows
\m{
\langle \Xi\rangle = \left(1+\sum_{j=1}^N |\xi_j|^2\right)^{1/2} \quad \mbox{for\, $\Xi = (\xi_1,\dots,\xi_N) \in \R^{nN}$ with $\xi_j \in \R^n$, \,j=\emph{1},\dots, N.}
}
The parameter $m$ is referred to as the \emph{order} or \emph{decay} of the amplitude.
\end{defn}
In what follows, we shall also use  \m{
\widehat{f}(\xi)=\int_{\R^n} f(x)e^{-i x.\xi}\, \dd x
}
as the definition of Fourier transform of $f$.
We now consider multilinear OIOs of the form 
\nm{eq:foidef}{
T^\Phase_\sigma(f_1,\dots,f_N)(x) = \int_{\R^{nN}} \sigma(x,\Xi)\prod_{j=1}^N \widehat{f}_j(\xi_j)\, e^{i\Phase(x,\Xi)}\, \ddd\Xi,
}
where $\sigma\in S^m(n,N)$ and $\ddd\Xi:= \dd\Xi/(2\pi)^{nN}.$

The main goal here is to show that the operator $T^\Phi_\sigma$, initially defined by \eqref{eq:foidef} for $f_1,\dots,f_N \in \mathscr{S}$ (the Schwartz class), extends to a bounded multilinear operator from $X^{p_1}\times \ldots \times X^{p_N}$ to $X^{p_0}$, where $X^{p_j}$ are certain Banach or quasi-Banach spaces. Now, in the case that $\frac{1}{p_0}=\sum_{j=1}^N\frac{1}{p_j},$ we shall refer to the corresponding regularity results as \emph{H\"older-type} (HT for short), otherwise \emph{non-H\"older-type} (NHT for short).\\

Given $\sigma \in S^{m}(n,N)$, the phases $\Phi$ in $T_a^{\Phi}$ for which regularity results are currently known take of the following forms:
\begin{enumerate}[label=\textit{\alph*})]
    \item $N=2$, $\Phase(x,\Xi) = \lambda\phase_0(\Xi)+ \sum_{j=1}^2 x\cdot\xi_j$,\, $\Xi  \in \R^{2n}$ and $\lambda$ a parameter;
    \item $N=2$, $\Phase(x,\Xi) =  \sum_{j=1}^{2}\varphi_j (x, \xi_j),$ \, $\Xi  \in \R^{2n}$; and
    \item $N\geq 1,$ $\Phase(x,\Xi) = \phase_0(\xi_1+\dots+\xi_N) + \sum_{j=1}^N (x\cdot\xi_j+\phase_j(\xi_j))$.
\end{enumerate}
For the phase functions of the form \textit a), one is aiming at  non-H\"older-type boundedness of $T_a^{\Phi}$ where part of the goal is also to obtain optimal powers of $\lambda$ in the boundedness estimates. In this case Bernicot and Germain \cite{BernicotGermain1} proved optimal global NHT regularity results in Lebesgue spaces, under suitable conditions on the rank of various Hessians of $\varphi_0$. Their analysis also accommodates quadratic phases.
 
For case \textit b) D. Rule, S. Rodr\'iguez-L\'opez and W. Staubach \cite{Monster} proved optimal HT local regularity results, under the conditions that the mixed Hessian of the phase functions $\varphi_j (x, \xi_j)$ are non-vanishing (the non-degeneracy condition) and that each of these phases are positively homogeneous of degree one in $\xi_j.$ Note that this case only accommodates examples that are relevant to the study of nonlinear wave equation. In \cite{RRS} it was shown that for bilinear operators where the phase functions are also allowed to behave quadratically, one can prove an $L^2 \times L^2 \to L^1$ boundedness result. A. Bergfeldt and W. Staubach \cite{BS} extended this to the case of globally defined multilinear operators and all possible Banach target spaces.

For case \textit c) Rule, Rodr\'iguez-L\'opez and Staubach \cite{RRS2} proved optimal HT global regularity results in the general multilinear case, under the condition that the phase functions $\varphi_j$ are positively homogeneous of degree one.  In this context, only the case of Banach target spaces were investigated.\\

\subsection{Synopsis of the results of the paper}
Given our previous discussions, there are quite a few problems that remain in the context of the regularity of oscillatory integral operators. Generally speaking, these problems are related to the nature of the amplitudes $a(x, \Xi)$ and that of the phase functions $\Phi(x, \Xi)$ for which one can prove various regularity results. For example one could lower the regularity of amplitudes or allow the phases to depend in a particular way on the spatial and/or frequency variables. In this paper we have chosen to look at the problem of global regularity for multilinear operators with phase functions of form $c)$ above, partly because of its relevance to the method of space-time resonance and partly because it is a tractable halfway house that should lead to an understanding of more general phases.

To implement our agenda, and motivated by examples related to dispersive PDEs, we consider the following class of phase functions:
\begin{defn}\label{defn:linear phase}
Let $0<s<\infty$. A function $\phase\colon \R^n \to \R$ which belongs to $\mathcal{C}^{\infty}(\mathbb{R}^n\setminus \{0\})$ and satisfies
\begin{equation}\label{eq:simplified phase}
\abs{\partial^\alpha\varphi(\xi)} \leq c_\alpha \abs\xi^{s-\abs\alpha} \text{ for }\, \xi\neq 0 \,\text{and}\, \abs\alpha \geq 0,
\end{equation} is called a \emph{phase function} \emph{(or phase)} of order $s$.
\end{defn}
We note that the case of the water wave equation corresponds to the case $s=\frac{1}{2}$, capillary waves to the case $s=\frac{3}{2}$, the Sch\"odinger equation to the case $s=2$ and the Airy equation to the case $s=3$.
 
We shall say that $0<p_j \leq \infty$ satisfy the H\"older condition if
\begin{equation}\label{eq:hoelder}
    \frac{1}{p_0}=\sum_{j=1}^N\frac{1}{p_j}.
\end{equation}
Now defining the functions spaces $X^p$ as
\begin{equation} \label{defn:xp}
    X^p:=\begin{cases}
        h^p & \mbox{if $p\leq 1$}\\
        L^p & \mbox{if $1<p<\infty$}\\
        \bmo & \mbox{if $p=\infty$},
    \end{cases}
\end{equation}
where $L^p$ is the usual Lebesgue space, $h^p$ is the local Hardy space defined in Definition \ref{def:Triebel} below, and $\bmo$ is the dual space of $h^1$, and considering phase functions of the form
\begin{equation}\label{eq:phaseform}
    \Phase(x,\Xi) = \phase_0(\xi_1+\dots+\xi_N) + \sum_{j=1}^N (x\cdot\xi_j+\phase_j(\xi_j)),
\end{equation}
we have the following HT boundedness result.
\begin{thm}\label{thm:main FIO}
For integers $n, N \geq2$, let the exponents $p_j \in (\frac{n}{n+1},\infty]$ \emph{($j=0,\dots, N$)} satisfy \eqref{eq:hoelder}.
Moreover let 
\begin{equation}\label{ineq:criticalexponent1}
	m\leq -(n-1)\brkt{\sum_{j=1}^N\abs{\frac{1}{p_j}-\frac{1}{2}}+\abs{\frac{1}{p_0}-\frac{1}{2}}}.
\end{equation}

If $\sigma \in S^{m}(n,N)$ and $\Phi$ is of the form \eqref{eq:phaseform} with each phase $\phase_j$ being smooth outside the origin and positively homogeneous of degree one, then the multilinear operator $T^\Phi_\sigma$ initially defined by \eqref{eq:foidef} for $f_1,\dots,f_N \in \mathscr{S}$ $($the Schwartz class$)$, extends to a bounded multilinear operator from $X^{p_1}\times \ldots \times X^{p_N}$ to $X^{p_0}$. Moreover, the same result holds in case each $\varphi_j(\xi)$ is equal to $\langle \xi\rangle$, which is an inhomogeneous phase related to the Klein--Gordon equation, for the range $p_j \in (0,\infty]$.
\end{thm}
For the range $p_0\in[1, \infty]$ (the Banach target-spaces), this theorem for the case of homogeneous of degree one phase functions (which is the case of the wave-equation) was proven in \cite{RRS2}.  Therefore Theorem \ref{thm:main FIO} extends our previous result to the quasi-Banach target-spaces as well as to the Klein--Gordon case. Note also that the admissible dimensions in the case of $N\geq 2$ are necessarily greater than or equal to two (see \cite{Monster}), however if $N=1$ then of course $n=1$ is also allowed, since this is just the well-known boundedness result for linear Fourier integral operators \cite{Miyachi}, \cite{Peral}. 
 Our second result HT boundedness result is the following.
\begin{thm}\label{thm:main}
For integers $N, n \geq1$, and a real number $s\in (0,\infty)$, assume that the exponents $p_j \in (\frac{n}{n+\min(1,s)},\infty]$ \emph{($j=0,\dots, N$)} satisfy \eqref{eq:hoelder}.
Suppose also that $\sigma \in S^{m}(n,N)$ and $\Phi$ is of the form \eqref{eq:phaseform} with each phase $\phase_j$ \emph{($j=0,1,\dots,N$)} of order $s$ and
\begin{equation}\label{ineq:criticalexponent2}
	m\leq -sn\brkt{\sum_{j=1}^N\abs{\frac{1}{p_j}-\frac{1}{2}}+\abs{\frac{1}{p_0}-\frac{1}{2}}}.
\end{equation}
Then the multilinear operator $T^\Phi_\sigma$ initially defined by \eqref{eq:foidef} for $f_1,\dots,f_N \in \mathscr{S}$, extends to a bounded multilinear operator from $X^{p_1}\times \ldots \times X^{p_N}$ to $X^{p_0}$.  Moreover, if the functions $\varphi_j$ are all in $\mathcal{C}^{\infty}(\mathbb{R}^n)$  $($the Schr\"odinger case is such an example$)$, then the ranges of the exponents $p_j$ in the theorem can be extended to $\in (0,\infty]$.
\end{thm}

\begin{rem}
If in Theorems \emph{\ref{thm:main FIO}} and \emph{\ref{thm:main}}, the phase function $\varphi_0=0$ $($and $n>1$ in the case of multilinear \emph{FIOs}$)$, then the order of the decay $m$ can be improved by just removing the term $-sn|1/p_{0} -1/2|$ $($or $-(n-1)|1/p_{0} -1/2|$$)$ from the $m$'s given in those theorems.
\end{rem}
Theorem \ref{thm:main} has no predecessor in the literature and covers the cases of water wave, capillary wave, Schr\"odinger, Korteweg--de\,Vries and many other higher order dispersive equations. Moreover, this result, in contrast to Theorem~\ref{thm:main FIO}, applies in all dimensions, when $N\geq 1$.

In proving Theorems \ref{thm:main FIO} and \ref{thm:main}, we make use of several global boundedness results: Those for linear Klein--Gordon equations, proved by J.~Peral \cite{Peral} (for $X^p$ with $1<p<\infty$); those for linear wave equations, proved by S.~Rodr\'iguez-L\'opez, D.~Rule and W.~Staubach \cite{RRS} (for $X^p$ with $n/(n+ 1)<p\leq \infty$); and those for higher order equations, proved by A.J.~Castro, A.~Israelsson, W. ~Staubach and M.~Yerlanov \cite{CISY} (for $X^p$ with $n/(n+\min(1,s))<p\leq \infty$).

The methods involved in proving the multilinear results in the realm of Banach spaces are essentially the same as the ones used by us to prove the boundedness of multilinear FIOs in \cite{RRS}, which are based on non-trivial extensions of the Coifman--Meyer methods in \cite{CM} to the case of multilinear operators with nonlinear phase functions.

{Thus, one writes the multilinear operator as a sum of operators whose amplitudes have specific support properties in the frequency variable $\Xi$. One term has compact frequency support, and for the other terms one has either that some $|\xi_j|$ dominates $\Xi,$ or that $\abs{\xi_j} \approx \abs{\xi_k}$ for certain $j$ and $k$ on the support of the amplitude in question.  Thereafter one identifies the end-points that are needed to apply complex interpolation and proving these end-point results creates in turn a number of cases, which are dealt with in accordance with whether the target spaces are Banach or quasi-Banach.

In the case of Banach target spaces and for the term that is compactly supported in the frequency, and those where $|\xi_j|$ dominates $\Xi,$  the machinery of \cite{RRS2} can be used without difficulty. However for the parts where $\abs{\xi_j} \approx \abs{\xi_k}$  (and for the target spaces $\bmo$ and $L^2$), one needs a result, provided in  Proposition \ref{lem:smallcarlnorm}, that demonstrates how certain oscillatory integral operators give rise to Carleson measures, with an estimate on their Carleson-norms. 
With this result at hand, the rest of the analysis is as in the case of multilinear FIOs in \cite{RRS2}.

The major hindrance to overcome here is that in the realm of quasi-Banach spaces, all Coifman--Meyer-type approaches, including the ones used in \cite{RRS2} or \cite{Monster} fail because of the impossibility of using duality arguments. Thus to prove results in the quasi-Banach realm, it behoves us to use a different method, and this is one of the novelties of the approach developed in this paper. To obtain the end-point results of this paper, our approach will be mainly based on various vector-valued inequalities. To our knowledge, using this type of estimates to derive estimates for multilinear oscillatory integral operators is new. The treatment that we describe here is rather technical, however it is fairly general in its nature and can be used in other contexts as well. We should mention, however, that this approach requires some degree of decay in the terms that represent the portion of multilinear operators where $\abs{\xi_j} \approx \abs{\xi_k}$. As such, the case of $L^{2}$-target spaces can not be subsumed in the quasi-Banach methods due to exactly that lack of decay. In addition to this, a lack of a convenient vector-valued characterisation for $\bmo$ means the method developed here also can not be applied in the case of a $\bmo$-target space. Fortunately though, the $L^{2}$ and $\bmo$-target space cases can be handled by the strategies mentioned earlier so that, in the end, we arrive at all the desired results for both Banach and quasi-Banach targets, albeit with a slightly longer proof than one might have hoped.}\\

The main motivation for our work was provided by a series of papers of F.~Bernicot and P.  Germain in \cite{BernicotGermain1,BernicotGermain2, BernicotGermain3} regarding coupled systems of dispersive PDEs, where the authors derived bilinear dispersive estimates in dimension 1, 2 and 3, for these systems in light of the method of space-time resonances. To briefly recall the setting of Bernicot-Germain's investigation,
let $\zeta(\Xi)$ be a smooth symbol
and let $T_ζ$ be the associated multilinear paraproduct  defined by 
\begin{equation}
\label{defn:tm}
T_\zeta( f_1, \dots,f_N)(x):=\int_{\R^{nN}}  \zeta(\Xi)\prod_{j=1}^N\left(\widehat{f}_j(\xi_j)e^{ix\cdot\xi_j}\right) \dd\Xi,
\end{equation}
where $\xi_j\in\R^n$ ($j=1,\dots,N$) and $\Xi=(\xi_1,\dots,\xi_N)\in \R^{nN}$. Furthermore, for $j=0,\dots,N$, let
\[
    \varphi_j(D)\, f(x)= \int_{\R^n}   \varphi_j(\xi)\,\widehat{f}(\xi)\,e^{ix\cdot\xi}\, \ddd\xi,
\]
where $\ddd \xi$ denotes the normalised Lebesgue measure ${\dd \xi}/{(2\pi)^n}$.
Consider now the coupled system of dispersive equations
\begin{equation*}
 \label{eq:dispersive}
\left\{ \begin{array}{l} i\partial_t u +  \varphi_0(D)\, u = T_{\zeta}\left( v_1,\dots, v_N\right)  \\
i\partial_t v_j +  \varphi_k(D)\, v_j = 0,\,\,\, j=1,\dots, N \\
\end{array} \right.
\quad \mbox{with} \quad
\left\{ \begin{array}{l} u(0,x) = 0  \\ v_j(0,x) = f_j (x), \,\,\, j=1,\dots, N. \end{array} \right.  
\end{equation*}
The functions $u$ and $v_k$ are complex valued, and each $f_k$ maps $\mathbb{R}^n$ to $\mathbb{C}$. %

{The above system is used in order to study the nonlinear interaction of free waves, as a first step towards understanding a nonlinear dispersive equation $i\partial_t u + \varphi(D) u=F(u)$, with a suitable nonlinearity. Thus given $f_j$ in some function spaces, one would like to understand the behaviour of $u$ in some other function spaces. 
}

Using this setting and our estimates for multilinear oscillatory integrals we are able to establish the validity of the following regularity theorem.

\begin{thm}\label{thm:PDE}
Let $s\in(0,\infty)$, $\sigma_k\geq 0$, $k=1, \dots, \, N$, $\varkappa=\min \sigma_k$, $p_j\in (1,\infty),$ $j=0, \dots, \, N$, satisfying the H\"older condition \eqref{eq:hoelder}, and assume that $\varphi_k \in \mathcal{C}^{\infty}(\R^n \setminus 0)$ are positively homogeneous of degree $s$, and $f_k\in H^{\sigma_k, p_k}$. Assume further that $T_\zeta$ is the multilinear multiplier given by \eqref{defn:tm} with symbol $\zeta(\Xi)\in S^{m_\zeta} (n, N)$ and set $m_c(s):=-ns\sum_{j=0}^N \abs{\frac 1 {p_j}-\frac 1 {2}}$ for $s\neq 1$, and $m_c(1)=-(n-1)\sum_{j=0}^N \abs{\frac 1 {p_j}-\frac 1 {2}}$. Then for any $q\in [1,\infty]$ and any $T>0$, there exists a constant $C_T>0$ such the solution $u(t, x)$ satisfies the regularity estimate
\[
   \Vert u \Vert_{L^{q}([0,T])\,H^{\varkappa+m_c(s)-m_{\zeta},p_0}(\R ^n)}\leq C_{T} \prod_{j=1}^{N}\Vert  f_j\Vert_{H^{\sigma_j, p_j}},
\]
provided that $\varkappa +m_c(s)-m_{\zeta}\geq 0$ $($which is needed in order to land in a space of functions rather than a space of distributions$)$.
\end{thm}
Here for $1<p<\infty$, $\sigma\in \R$, $H^{\sigma,p}=\{ f\in \mathscr{S}'; \, (1-\Delta)^{s/2}f\in L^p(\R^n)\}$ is the $L^p$-based Sobolev space with the norm
$\Vert f\Vert_{H^{\sigma, p}}:= \Vert (1-\Delta)^{s/2}f\Vert_{L^p}$.\\

The paper is organised as follows. In Section \ref{subsection:definitions} we recall the basic notions and tools from Fourier analysis and state some fairly general results that will also be used in the proof of Theorems \ref{thm:main FIO} and \ref{thm:main}. In Section \ref{sharpness section} we briefly discuss the sharpness of the order of the decay of the operators in the bilinear setting. In Section \ref{section main linear estim} we state and prove several results in the vector-valued setting for linear OIOs as well as a key proposition regarding the OIOs giving rise to Carleson measures. Section \ref{sec:freqdecom} recalls briefly the frequency decomposition that was introduced in \cite{RRS2}, and which will be used throughout the paper. In Section \ref{endpoint cases} we briefly discuss the endpoint cases that are going to be considered in the Banach-target case. Section \ref{sec main multilinear} contains the proofs of Theorems \ref{thm:main FIO} and \ref{thm:main}. Finally Section \ref{space-time} is devoted to the proof of Theorem \ref{thm:PDE} on the Sobolev regularity of the solutions to coupled systems of dispersive partial differential equations.

\section{Definitions and Preliminaries}\label{subsection:definitions}
Here we collect all the definitions and basic results that will be used in the forthcoming sections, in order to make the paper essentially self-contained. 

We shall denote constants which can be determined by known parameters in a given situation, but whose values are not crucial to the problem at hand, by $C$ or $c$, sometimes adding a subscript, for example $c_\alpha$, to emphasis a dependency on a given parameter $\alpha$. Such parameters are those which determine function spaces, such as $p$ or $m$ for example, the dimension $n$ of the underlying Euclidean space, and the constants connected to the seminorms of various amplitudes or phase functions. The value of the constants may differ from line to line, but in each instance could be estimated if necessary. We also write $a\lesssim b$ as shorthand for $a\leq Cb$ and $a\approx b$ when $a\lesssim b$ and $b\lesssim a$. By
\m{
B(x,r) := \{ y\in\R^n \, : \, |y-x| < r\}
}
we denote the open ball of radius $r > 0$ centred at $x\in\R^n$.

We also recall the definition of the \emph{Littlewood--Paley} partition of unity which is a basic tool in harmonic analysis and theory of partial differential equations.

\begin{defn}\label{def:LP}
Let $\vth \colon \R^n \to \R$ be a positive, radial, radially decreasing, smooth cut-off function which satisfies $\vth(\xi) = 1$ if $|\xi| \leq 1$ and $\vth(\xi) = 0$ if $|\xi| \geq 2$. We set $\vartheta_0:=\vartheta$ and 
$$
    \vartheta_j(\xi) := \vartheta \left (2^{-j}\xi \right )-\vartheta(2^{-(j-1)}\xi),
$$
for integers $j\geq 1$. Then one has the following \emph{Littlewood--Paley partition of unity}:
\[\label{eq:littlewoodpaley}
   \sum_{j=0}^\infty \vartheta_j(\xi) = 1 \quad \text{\emph{for all }}\xi\in\R^n.
\]
\end{defn}
Using the definition above, let $s \in {\R}$ and $0< p <\infty$, $0< q \leq\infty$. The Triebel--Lizorkin space is defined as
	\[\label{TLspace}
	   F^s_{p,q}(\R^n)
	:=
	\Big\{
	f \in {\mathscr{S}'}(\R^n) \,:\,
	\|f\|_{F^s_{p,q}(\R^n)}
	:=
	\Big\|\Big\{\sum_{j=0}^{\infty} 2^{jqs}\left|\vartheta_j(D) f\right|^{q}\Big\}^{1 / q}\Big\|_{L^{p}(\R^n)}<\infty
	\Big\}, 
	\]
where $\mathscr{S}'(\R^n)$ denotes the space of tempered distributions.	

In our analysis of the boundedness of oscillatory integral operators which is based on multilinear interpolation, the end-points often involve local Hardy spaces which were introduced by D. Goldberg \cite{Gold}. One of the main advantages of these spaces is that they are mapped into themselves under the action of the linear oscillatory integral operators that are considered in this paper.
\begin{defn}\label{def:Triebel}
The \emph{local Hardy space} $h^p(\R^n)$, $($$0<p<\infty$$)$ is the Triebel--Lizorkin space $F^{0}_{p,2}$ $($see, for example \cites{Triebel}$)$ with the norm

\begin{equation}\label{norm in loc hardy}
    \|f\|_{{h}^p(\R^n)} \approx \norm{\vartheta_{0}(D)f}_{L^p(\R^n)}+\norm{\left(	\sum_{j=1}^\infty|\vartheta_j(D)f|^{2} \right)^{\frac{1}{2}}}_{L^p(\R^n)}.
    \end{equation}
%

%
\end{defn}

Note that the usual Hardy space $\mathscr{H}^p(\R^n)$ is defined the condition
\begin{equation*}\label{eq:hp basic}
	\|f\|_{{\mathscr{H}}^p} := \brkt{\int \sup_{t>0} \abs{\vartheta(tD) f(x)}^p\dd x}^{\frac{1}{p}} <\infty.
\end{equation*}
The dual of $\mathscr{H}^1$ is the John--Nirenberg space of functions of bounded mean oscillation $\mathrm{BMO}$, which consists of
all functions $f\in L^1_{\mathrm{loc}}$ 
such that
\[\label{defn:bmonorm}
 \|f\|_{\mathrm{BMO}}
:=\sup_{Q}\frac{1}{|Q|}
\int_{Q}|f(x)-\mathrm{avg}_Q f|\, dx
<\infty,   
\]
where $\mathrm{avg}_Q f=|Q|^{-1}\int_Q f$, and $Q$ ranges over cubes in $\R^n$. The dual of the local Hardy space $h^1$ is the \textit{local} $\BMO$ space, which is denoted by $\bmo$ and consists of locally integrable functions that verify
\begin{equation*}\label{defn:bmo}
\Vert f\Vert_{\bmo}\approx \Vert f\Vert_{\BMO}+ \Vert\vartheta (D) f\Vert_{L^{\infty}}<\infty,
\end{equation*}
where $\vartheta$ is the cut-off function introduced in Definition \ref{def:LP}. 

In the analysis of multilinear operators, a basic tool is a certain type of measure whose definition we now recall. A Borel measure $\dddd\mu(x,t)$ on $\R^{n+1}_+$ is called a \emph{Carleson measure} if
\m{
\Vert \dddd \mu\Vert_{\mathcal{C}} := \sup_Q \frac{1}{|Q|}\int_0^{\ell(Q)}\int_Q |\dddd\mu(x,t)| < \infty
}
where the supremum is taken over cubes $Q \subset \R^n$ and $\ell(Q)$ denotes the side length of $Q$ and $|Q|$ its Lebesgue measure. The quantity $\Vert \dddd \mu\Vert_{\mathcal{C}}$ is called the \emph{Carleson norm} of $\dd \mu$. An equivalent norm is given if cubes are replaced with balls. In this paper we are exclusively interested in Carleson measures which are supported on lines parallel to the boundary of $\R^{n+1}_+$. More precisely, in what follows all Carleson measures will be supported on the set
\m{
E := \{(x,t) \, : \, \mbox{$x\in\R^n$ and $t=2^{-k}$ for some $k\in\Z$}\}
}
so they take the form
\m{
\sum_{k\in\Z}\dddd\mu(x,t)\delta_{2^{-k}}(t),
}
where $\delta_{2^{-k}}(t)$ is a Dirac measure at $2^{-k}$. This will be assumed throughout without further comment.

The following basic results concerning the Carleson measure and the quadratic estimate are very useful in the context of multilinear operators. See E. M. Stein \cite{S} for the proofs.
{{
\begin{lemma}\label{Coifman-Meyerslem}
Let $\dd \mu(x, t)$ be a Carleson measure. Then 
if $\varphi$ satisfies $|\varphi(x)|\lesssim \langle x\rangle^{-n-\varepsilon}$ \emph{(for some $0<\varepsilon<\infty$), then}
\begin{equation}\label{ineq:carll2}
    \sum_{k}\int _{\R^n}|\varphi(2^{-k}D) f(x)|^2 \, \dd \mu(x, 2^{-k})\leq C_n \norm{\dddd \mu}_{\mathcal{C}} \norm{f}_{L^2}^2,
\end{equation}
and if $\varphi$ is a bump function supported in a ball near the origin with $\varphi(0)=1$ then one also has
\begin{equation}\label{ineq:carlh1}
    \sum_{k}\int _{\R^n}|\varphi(2^{-k}D) f(x)| \, \dd \mu(x, 2^{-k})\leq C_n \norm{\dd \mu}_{\mathcal{C}} \norm{f}_{h^1}.
\end{equation}
If $\varphi \in {\mathscr{S}}$ is such that $\varphi(0) = 0$, then
\begin{equation}\label{ineq:quadraticestimate}
\sum_k \int \abs{\varphi(2^{-k}D)f(x)}^2 \dd x \lesssim \norm{f}_{L^2}^2.
\end{equation}

\end{lemma}}}

In our investigations we will also confront three types of maximal operators. The first one is the  Hardy--Littlewood maximal operator

 \begin{equation*}\label{HLmax}
  \mathcal{M} f(x):=\sup _{B\ni x}\frac{1}{|B|} \int_{B}|f(y)| d y,   
 \end{equation*} where the supremum is taken over all balls $B$ containing $x$. For \(0<p<\infty\), one also defines 
\(\mathcal{M}_{p} f(x):=\left(\mathcal{M}\left(|f|^{p}\right)\right)^{1 / p}\).

The second one is J. Peetre's maximal operator \cite{Triebel}.
\begin{equation}\label{Peetremax}
  \mathfrak{M}_{a,b}(f)(x):=\Big\Vert  \frac{{f(x-\cdot)}}{\brkt{1+b\abs{\cdot}}^a}\Big\Vert_{L^{\infty} } 
\end{equation}
where $0<a, b<\infty$. For any $x\in \R^n$, $f\in\mathscr{S}'$ with $\supp \hat{f}\subset \{\xi;\,|\xi| \leq 2b \}$ and $a\geq \frac{n}{p}$ one has that

\begin{equation}\label{hl bounds peetre}
\mathfrak{M}_{a, b} u(x) \lesssim \mathcal{M}_{p} u(x).
\end{equation}

The third type of maximal operator that will be used in this paper is B.J. Park's maximal operator \cite{Park}: For $j \in \mathbb{Z}$, $s>0$ and $0<p \leq \infty$
\begin{equation}\label{Parkmax}
\mathfrak{M}_{s, 2^{j}}^{p} f(x):=2^{j n / p}\left\|\frac{f(x-\cdot)}{\left(1+2^{j}|\cdot|\right)^{s}}\right\|_{L^{p}}.
\end{equation}
Park's maximal operator has the following properties: If \(0<p<\infty\) and \(s>n / p\), then
\begin{equation}
    \mathfrak{M}_{s, 2^{j}}^{p} f(x) \lesssim \mathcal{M}_{p} f(x),
\end{equation}
uniformly in \(j \in \mathbb{Z}.\)
Moreover if the set of all dyadic cubes in $\R^n$
is denoted by $\mathcal{D}$, and for each $j\in \Z$ one denotes the elements of $\mathcal{D}$ with side length $2^{-j}$ by $\mathcal{D}_j$, then for every dyadic cube $J \in \mathcal{D}_j$ and for every $s>0$, $0<p<\infty$ and $f$,
\begin{equation}\label{Parks inequality}
    \sup_{y\in J} \mathfrak{M}_{s,2^j}^p f(y)\lesssim \inf_{y\in J} \mathfrak{M}_{s,2^j}^p f(y),
\end{equation}
with constants independent of $f$ and $j$.

Using the maximal operator $\mathfrak{M}_{a,b}$, Park has given a useful characterisation of the Hardy and BMO spaces, in the following theorem.

\begin{thm} \cite{Park}.\label{Parks thm}
Let $\Lambda\in \mathscr{S}$ be a function whose Fourier transform is supported in the annulus \(1 / 2\leq |\xi|\leq 2\) and set \(\widehat{\Lambda}\left(\cdot / 2^{j}\right)=\widehat{\Lambda_{j}} \) so that one has the partition of unity
\(\sum_{j \in \mathbb{Z}} \widehat{\Lambda}_j\left(\xi\right)=1\) for \(\xi \neq 0 .\)  Assume that \(0<p \leq \infty, 0<q \leq \infty,\, 0<\gamma<1\), and \(s>n / \min (p, 2, q) .\) Then
for each dyadic cubes \(Q \in \mathcal{D}\), there exists a proper measurable subset \(S_{Q}\) of \(Q\), depending
on \(\gamma, s, q\) and $f$, such that \(\left|S_{Q}\right|>\gamma|Q|\) and
\begin{equation*}\label{Parkestim}
    \|f\|_{Y^{p}} \approx\left\|\left\{\sum_{Q \in \mathcal{D}_{j}}\left(\inf _{y \in Q} \mathfrak{M}_{s, 2^{j}}^{q}\left(\Lambda_{j} * f\right)(y)\right) \chi_{S_{Q}}\right\}_{j \in \mathbb{Z}}\right\|_{L^{p}\left(\ell^{2}\right)}
\end{equation*}
where \(Y^{p}=\mathscr{H}^{p}\) for \(0<p<\infty\) and \(Y^{\infty}=\mathrm{B M O}.\)
\end{thm}

Now in connection to the Hardy--Littlewood maximal operator defined above, a useful device in proving multilinear estimates is the Fefferman--Stein vector-valued maximal inequality \cite[Theorem 1]{FS2}, which states that for $r<p,$ $q<\infty$, or \(0<p<\infty,\, q=\infty\) or for \(p=q=\infty\), one has
\begin{equation}\label{FSvector}
\left\|\left\{\mathcal{M}_{r} f_{j}\right\}_{j \in \mathbb{Z}}\right\|_{L^{p}(\ell^q)} \lesssim\left\|\left\{f_{j}\right\}_{j \in \mathbb{Z}}\right\|_{L^{p}(\ell^q)}. 
\end{equation}
The following theorem gives a corresponding vector-valued inequality involving Park's maximal operator.
\begin{thm}\label{Park2}\cite{Park}.
Let \(0<p, q, r \leq \infty\) and \(s>n / \min (p, q, r) .\) Suppose that the Fourier
transform of \(f_{j}\) is supported in a ball of radius \(A 2^{j}\) for some \(A>0 .\)
Then for \(0<p<\infty\) and \(0<q \leq \infty\) or for \(p=q=\infty\), one has
\begin{equation}
\left\|\left\{\mathfrak{M}_{s, 2^{j}}^{r} f_{j}\right\}_{j \in \mathbb{Z}}\right\|_{L^{p}(\ell ^q)} \lesssim\left\|\left\{f_{j}\right\}_{j \in \mathbb{Z}}\right\|_{L^{p}(\ell^q)}
\end{equation}
\end{thm}
We will also need the following vector valued inequality due to H. Triebel \cite[Theorem 2, Section 2.4.9]{Triebel}.
\begin{thm}\label{triebels thm}
If $G_k$ is a sequence of functions with $\mathrm{supp}\,\widehat{G_k}\subset B(0, 2^{k} R)$, for $k=0, 1,\dots$ and $R\geq 1$, then  for $0<r<\infty$ and $0< q<\infty$ one has the following vector-valued inequality$:$ For $\mathfrak{m}\in H^{\alpha}(\R^n)$ $($the Sobolev space $H^{\alpha,2}$ of order $\alpha$ defined in the introduction section$)$, and $\widehat{\mathfrak{m}(2^{-k}D)f}(\xi)= \mathfrak{m}(2^{-k}\xi)\hat{f}(\xi)$, with 
\[
    \alpha>n\left(\frac{1}{\min (1, r, q)}-\frac{1}{2}\right),
\] 
there is a constant $C>0$ independent of $R$ and $G_k$'s, such that
\begin{equation}\label{t_Triebels vector valued}
\norm{\Big\{\mathfrak{m}(2^{-k}D)
 G_k\Big\}_{k\in \Z}}_{L^r (\ell^q)}\leq C\Vert \mathfrak{m}\Vert_{H^{\alpha}} \Big\Vert\Big\{G_k\Big\}_{k\in \Z}\Big\Vert_{L^r(\ell^q)}.
 \end{equation}

\end{thm}
We note that the multilinear amplitudes defined in Definition \ref{defn:amplitudes} reduce to the classical H\"ormander classes $S^m$ of \emph{amplitudes} (or \emph{symbols}) in the case $N=1$, that is to say $S^m = S^m(n,1)$. The linear OIOs are the special case of \eqref{eq:foidef} when $N=1$, in which case we have
\begin{equation}\label{defn:linear OIO}
T_a^{\phase}f(x) := \int_{\R^n}e^{ix\cdot\xi +i\varphi(\xi)}a(x,\xi)\widehat{f}(\xi) \ddd \xi,
\end{equation}
for a given amplitude $a \in S^m$ and phase function $\phase$. In the proofs in the forthcoming sections we will also use the notion of \emph{multilinear pseudodifferential operators} which are operators of the form
\[\label{eq:mpseudodef}
T_\sigma(f_1,\dots,f_N)(x) = \int_{\R^{nN}} \sigma(x,\Xi)\prod_{j=1}^N \widehat{f}_j(\xi_j)\, e^{i\sum_{j=1}^N x\cdot\xi_j}\, \ddd\Xi.
\]

For the analysis of the low frequency portion of the operators, where ususally the singularity of the phase functions lie, we recall a linear result proved in \cite{CISY}, which established the $h^p$-boundedness of low-frequency portions of oscillatory integral operators, whose multilinear generalisations are considered here in this paper.

\begin{lem}\label{low freq lemma llf 1}
Let $s>0$,  $s_c:=\min(s,1)$, $a(x,\xi)
$ be a symbol that is compactly supported and smooth outside the origin in the $\xi$-variable and $\varphi(\xi)\in \mathcal{C}^{\infty}( \R^n \setminus \{0\}) $ be a phase function. Also assume that the following conditions hold:
\begin{align*}
&\begin{cases}\Vert \partial_{\xi}^{\alpha} a(\cdot,\xi)\Vert_{L^\infty(\R^n)}
\leq c_\alpha , 
&|\alpha|\geq 0,\\
|\partial^{\alpha}_{\xi} \varphi(\xi) |\leq c_{\alpha} |\xi|^{s-|\alpha|}, & \abs\alpha\geq 0,
\end{cases}\label{eq:lf2}
\end{align*}
{for $\xi\neq0$} and on the support of $a(x,\xi).$ Let 
\[
    K(x,y) := \int_{\R^n} a(x,\xi)\, e^{i\varphi(\xi){\color{blue}\boldsymbol+}i(x-y)\cdot \xi}\ddd\xi.
\]
Then one has:
\begin{enumerate}
    \item\label{lem2.6:parti} $\displaystyle \abs{K(x,y)}\lesssim \langle x-y\rangle ^{-n-\eps s_c}$ for any $0\leq \eps <1$.
    \item For every $r \in ( n/(n+\eps s_c),1]$ one has, for every $f\in \mathscr{S}'$ with frequency support inside the unit ball and $T_{a}^{\varphi}$ defined as in \eqref{defn:linear OIO}, that
\begin{equation*}
     \abs{T_{a}^{\varphi} f(x)}\lesssim \mathcal{M}_r f (x), \quad x \in \mathbb{R}^n.
\end{equation*}

\item For every $\frac{n}{n+s_c}<p\leq \infty$, and all $f\in X^p$, 
\[
    \norm{T_{a}^\varphi f}_{X^p}\lesssim \norm{f}_{X^p}.
\]
\end{enumerate}
\end{lem}
\begin{proof} The proof of the first statement can be found in \cite{CISY}*{Lemma 4.3}.

For the second statement we can apply \eqref{lem2.6:parti} to obtain that
 \begin{equation*}\label{poitwise estimate for the lowfrequency part}
 |T_{a}^{\phase}f(x)|\lesssim |(\vartheta (D) f) \ast \langle \cdot \rangle^{-n-\varepsilon s_c}|\lesssim \mathcal{M}_r (\vartheta (D) f)(x)
 \end{equation*}
 for all $f\in \mathscr{S}$, $r\in (\frac{n}{n+\varepsilon s_c},1]$ and $\varepsilon \in (0,1)$.
 
We can prove the third statement by choosing $\frac{n}{n+s_c}<r<p$ and making use of the boundedness of the maximal operator $\mathcal{M}$ on $L^{p/r}$ to obtain
\m{
\norm{T_{a}^{\phase}f}_{h^p}\lesssim \Vert T^{\varphi}_{a} f\Vert_{L^p} \lesssim \Vert \mathcal{M}(|\vartheta (D) f|^r)\Vert_{L^{p/r}}^{1/r} \lesssim \Vert \vartheta(D) f\Vert_{L^{p}} \lesssim \Vert f\Vert_{h^p},
}
where the last inequality follows by \eqref{norm in loc hardy} in Definition \ref{def:Triebel}. 
In the case of $p=\infty$ for which $X^p=\mathrm{bmo}$, we just observe that the integral kernel of the adjoint of $T_{a}^{\phase}$ is given by
 $ \int_{\R^n} a(y,\xi)\, e^{-i\varphi(\xi)-i(x-y)\cdot \xi}\ddd\xi,$ for which one can deduce a similar decay estimate as in ($i$). Therefore by the same reasoning as above one has that $\norm{(T_{a}^{\phase})^{\ast}f}_{h^1}\lesssim \norm{f}_{h^1}$ and hence $T_{a}^{\phase}$ is bounded on $\bmo$.
\end{proof}

As was mentioned earlier, the proofs of Theorems \ref{thm:main FIO} and \ref{thm:main} also use the following linear results: 

\begin{thm}\label{linearhpthmfio}
Let  $m=-(n-1)\abs{\frac{1}{p}-\frac{1}{2}}$ and $\frac{n}{n+1}<p\leq\infty$. Then any \emph{FIO} of the form 
$$T_\sigma^{\phase}f(x)= \int_{\R^n} \sigma(x,\xi)\, e^{ix\cdot \xi +i\phase(\xi)} \widehat{f}(\xi) \ddd\xi ,$$
with an amplitude $\sigma(x,\xi)\in S^{m}$ and a real-valued phase function $\phase\in \mathcal{C}^{\infty}(\mathbb{R}^n\setminus \{0\})$ that is positively homogeneous of degree one, satisfies the estimate
\[
	\norm{T_\sigma^{\phase} f}_{X^p}\leq C\norm{f}_{X^p},
\]
where $X^p$ is defined in \eqref{defn:xp}. Moreover, the same result also holds for $0<p<\infty,$ if $\varphi(\xi)$ is equal to the inhomogeneous phase function $\langle \xi\rangle$ $($the case of the Klein--Gordon equation$).$
\end{thm}
\begin{proof}
For homogeneous phase functions, this result was established in \cite[Theorem 3.1]{RRS2}. For the proof for $\varphi(\xi)=\langle \xi \rangle$ we sketch an argument from \cite{IMS}. One first separates the amplitude $\sigma(x, \xi)$ into low and high frequency portions. For the low frequency part we have the result thanks to Lemma \ref{low freq lemma llf 1}, and for the high frequency part one can write $\sigma(x,\xi) e^{i\langle \xi \rangle}= \tilde{\sigma}(x,\xi) e^{i|\xi|}$ with $\tilde{\sigma}\in S^m$ and thereafter apply Theorem 3.1 from \cites{RRS2} once again.
\end{proof}
For other classes of OIOs, the following theorem was proven in \cite{CISY}, Theorem 3.5.
\begin{thm}\label{linearhpthmoio}
Let $0<s<\infty$, $m=-ns\abs{\frac{1}{p}-\frac{1}{2}}$ and $\frac{n}{n+\min(s,1)}<p\leq\infty$. Then any linear oscillatory integral operator
$$T_\sigma^{\phase}f(x)= \int_{\R^n} \sigma(x,\xi)\, e^{ix\cdot \xi +i\phase(\xi)} \widehat{f}(\xi) \ddd\xi ,$$
with an amplitude $\sigma(x,\xi)\in S^{m}$ and a phase function $\varphi$ satisfying \eqref{eq:simplified phase}, satisfies the estimate
\[
	\norm{T_\sigma^{\phase} f}_{X^p}\leq C\norm{f}_{X^p}.
\]
Moreover, if the phase function $\varphi$ is in $\mathcal{C}^{\infty}(\mathbb{R}^n),$ then the range of $p$ in the theorem can be extended to $\in (0,\infty]$.
\end{thm}

\section{On the sharpness of the orders of the operators}\label{sharpness section}

Here, building on the example in \cite{RRS2} and the sharpness results in \cite {Miyachi}, we construct examples which show the sharpness of  \cite[Theorem~2.7]{Monster} for certain values of the function space exponents. They also serve as examples which show the sharpness of our main results here (Theorems~\ref{thm:main FIO} and \ref{thm:main}) when the target space is $L^2$. As such, we consider the case of bilinear operators with $\varphi_0=0$, and the failure of $L^p \times L^q \to L^r$ boundedness (in the cases $p, q\leq 2$ and $p,q\geq 2$). At the very end of the section we consider the case of $\varphi_0\neq 0$ but only for $r=2$.

So let us first consider the operator 
\begin{equation*}
B(f,g)(x)= \int_{\mathbb{R}^{2n} } a(\xi,\eta)\widehat{f}(\xi) \, \widehat{g}(\eta)\, e^{ix\cdot(\xi+\eta)}\, e^{ i \phase(\xi) - i \phase(\eta)} \, \dd\xi \, \dd\eta,
\end{equation*}
with $\phase(\xi) = |\xi|^s$,
\begin{equation*}
    a(\xi,\eta) = \sum_{k=0}^\infty \vartheta_k(\xi)\overline{\vartheta_k(-\eta)}b_1(\xi)\overline{b_2(-\eta)},
\end{equation*}
and $b_j(\xi) = (1-\vartheta_0(\xi))|\xi|^{m_j}$ ($j=1,2$), so that $a \in S^{m}_{1,0}(n,2)$, with $m = m_1+m_2$.

The parameter $m$ and the order $s$ of $\phase$ will be specified later, but we have in mind that $m$ should fail to satisfy either \eqref{ineq:criticalexponent1} or alternatively \eqref{ineq:criticalexponent2} depending on $s$. We compute
\begin{equation}
\begin{aligned}\label{eq:sharpness1}
& B(f,\overline{g})(x) \\
&= \int_{\R^{2n} } \left(\sum_{k=0}^\infty \vartheta_k(\xi)\overline{\vartheta_k(-\eta)}b_1(\xi)\overline{b_2(-\eta)}\right)\widehat{f}(\xi) \, \widehat{g}(-\eta)\, e^{ix\cdot(\xi+\eta)}\, e^{ i \phase(\xi) - i \phase(\eta)} \, \dd\xi \, \dd\eta \\
&= \sum_{k=0}^\infty \left(\int_{\R^{n} } \vartheta_k(\xi)b_1(\xi)\widehat{f}(\xi) \, e^{ix\cdot\xi}\, e^{ i \phase(\xi)} \, \dd\xi\right) \left(\int_{\R^{n} } \overline{\vartheta_k(-\eta)b_2(-\eta)\widehat{g}(-\eta)} \, e^{ix\cdot\eta}\, e^{ -i \phase(\eta)} \, \dd\eta\right) \\
&= \sum_{k=0}^\infty \left(\int_{\R^{n} } \vartheta_k(\xi)b_1(\xi)\widehat{f}(\xi) \, e^{ix\cdot\xi}\, e^{ i \phase(\xi)} \, \dd\xi\right) \overline{\left(\int_{\R^{n} }\vartheta_k(\xi)b_2(\xi)\widehat{g}(\xi) \, e^{ix\cdot\xi}\, e^{ i \phase(\xi)} \, \dd\xi\right)}.
\end{aligned}
\end{equation}

\subsection{Fourier integral operators}
Consider $s=1$ and
\begin{equation*}
m = -(n-1)\brkt{\abs{\frac{1}{p}-\frac{1}{2}} + \abs{\frac{1}{q}-\frac{1}{2}}} + \varepsilon
\end{equation*}
for some $\varepsilon > 0$.

If $p,q \geq 2$ (so $2r\geq 1$) we choose
\begin{align*}
    \lambda_1 &= \frac{n+1}{2} - \frac{1}{p} + \frac{\varepsilon}{4}, \\
    \lambda_2 &= \frac{n+1}{2} - \frac{1}{q} + \frac{\varepsilon}{4}, \\
    m_1 &= -\frac{n-1}{2} + \frac{n}{2r} - \frac{1}{p} + \frac{\varepsilon}{2}, \quad \mbox{and} \\
    m_2 &= -\frac{n-1}{2} + \frac{n}{2r} - \frac{1}{q} + \frac{\varepsilon}{2}. 
\end{align*}
We see directly that $m = m_1 + m_2$ and if we define $\widehat{f}(\xi)= (1-\vartheta_0(\xi)) |\xi|^{-\lambda_1}e^{-i|\xi|}$ and $\widehat{g}(\xi)= (1-\vartheta_0(\xi)) |\xi|^{-\lambda_2}e^{-i|\xi|}$, fact (II-i) from \cite[page 302]{Miyachi} shows us that $f\in L^p$ and $g \in L^q$. We see also that
\begin{equation*}b_1(\xi)\widehat{f}(\xi)e^{ i \phase(\xi)}
 = b_2(\xi)\widehat{g}(\xi) e^{ i \phase(\xi)} = (1-\vartheta_0(\xi))^2|\xi|^{-n(1-1/(2r))+\varepsilon/4} =: \widehat{F}(\xi).
\end{equation*}
so we can compute from \eqref{eq:sharpness1} that
\begin{equation}\label{eq:sharpness2}
    B(f,\overline{g})(x) = \sum_{k=0}^\infty \left|\int_{\R^{n} } \vartheta_k(\xi)\widehat{F}(\xi) \, e^{ix\cdot\xi} \, \dd\xi\right|^2
= \sum_{k=0}^\infty \left|\vartheta_k(D)(F)(x)\right|^2.
\end{equation}
If we assume $B$ is bounded from $L^p\times L^q$ to $L^{r}$, then the Littlewood-Paley characterisation of $h^{2r}$ and the fact that $F$ is high-frequency localised, yield
\begin{equation*} \label{eq:sharpness}
\begin{aligned}
   \|F \|_{H^{2r}} \sim    \|F \|_{h^{2r}} &\lesssim \left\Vert\left(\sum_{k=1}^\infty \left|\vartheta_k(D)(F)\right|^2\right)^{1/2}\right\Vert_{L^{2r}} 
    = \left\Vert T(f,\overline{g}) \right\Vert_{L^{r}}^{1/2} \lesssim \|f\|_{L^p}^{1/2}\|g\|_{L^q}^{1/2}.
\end{aligned}
\end{equation*}
However, fact (II-i) from \cite[page 302]{Miyachi} shows us that $F \not\in H^{2r}$. 

So we arrive at a contradiction, and $B$ cannot be a bounded operator from $L^p\times L^q$ to $L^{r}$.

If $p,q \leq 2$ we can apply a similar argument but choose instead
\begin{align*}
    \lambda_1 &= n\brkt{1 - \frac{1}{p}} + \frac{\varepsilon}{4}, \\
    \lambda_2 &= n\brkt{1 - \frac{1}{q}} + \frac{\varepsilon}{4}, \\
    m_1 &= \frac{n-1}{2} - \frac{n}{p} + \frac{1}{2r} + \frac{\varepsilon}{2}, \quad \mbox{and} \\
    m_2 &= \frac{n-1}{2} - \frac{n}{q} + \frac{1}{2r} + \frac{\varepsilon}2.
\end{align*}
We still have that $m = m_1 + m_2$ and if we this time define $\widehat{f}(\xi)= (1-\vartheta_0(\xi)) |\xi|^{-\lambda_1}$ and $\widehat{g}(\xi)= (1-\vartheta_0(\xi)) |\xi|^{-\lambda_2}$, fact (II-ii) from \cite[page 302]{Miyachi} shows us that $f\in L^p$ and $g \in L^q$. We once again obtain \eqref{eq:sharpness2} but with
\begin{equation*}
     \widehat{F}(\xi) = (1-\vartheta_0(\xi))^2|\xi|^{-(n+1)/2 + 1/(2r) + \varepsilon/4}e^{-i|\xi|},
\end{equation*}
so the proof of fact (II-ii) from \cite[page 302]{Miyachi} reveals that $F(x) \sim (1-|x|)^{-1/(2r) - \varepsilon/4}$ as $|x| \to 1$, so again $F\not\in H^{2r}$. We have therefore shown, even for $p,q \leq 2$, $B$ is not a bounded operator from $L^p\times L^p$ to $L^{p/2}$.

\subsection{Oscillatory integral operators}
We consider now either $0 < s < 1$ or $s > 1$ and
\begin{equation*}
m = -sn\brkt{\abs{\frac{1}{p}-\frac{1}{2}} + \abs{\frac{1}{q}-\frac{1}{2}}} + \varepsilon
\end{equation*}
for some $\varepsilon > 0$.

If $p,q \geq 2$ we choose
\begin{align*}
    \lambda_1 &= n\brkt{1 - \frac{s}{2}} - n\frac{(1-s)}{p} + \frac{\varepsilon}{4}, \\
    \lambda_2 &= n\brkt{1 - \frac{s}{2}} - n\frac{(1-s)}{q} + \frac{\varepsilon}{4}, \\
    m_1 &= -sn\brkt{\frac{1}{2} - \frac{1}{p}} - n\brkt{\frac{1}{p} - \frac{1}{2r}} + \frac{\varepsilon}{2}, \quad \mbox{and} \\
    m_1 &= -sn\brkt{\frac{1}{2} - \frac{1}{p}} - n\brkt{\frac{1}{p} - \frac{1}{2r}} + \frac{\varepsilon}{2}. 
\end{align*}
then we can carry out an analogous argument to that above for FIOs with $\widehat{f}(\xi)= (1-\vartheta_0(\xi)) |\xi|^{-\lambda_1}e^{-i|\xi|^a}$ and $\widehat{g}(\xi)= (1-\vartheta_0(\xi)) |\xi|^{-\lambda_2}e^{-i|\xi|^a}$. We use (I-i) instead of (II-i) from \cite{Miyachi} to conclude that $f \in L^p$ and $g \in L^q$ but $B(f,g) \not\in L^r$.

If $p,q \leq 2$ we choose
\begin{align*}
    \lambda_1 &= n\brkt{1 - \frac{1}{p}} + \frac{\varepsilon}{4}, \\
    \lambda_2 &= n\brkt{1 - \frac{1}{q}} + \frac{\varepsilon}{4}, \\
    m_1 &= -sn\brkt{\frac{1}{2r} - \frac{1}{2}} - n\brkt{\frac{1}{p} - \frac{1}{2r}} + \frac{\varepsilon}{2}, \quad \mbox{and} \\
    m_1 &= -sn\brkt{\frac{1}{2r} - \frac{1}{2}} - n\brkt{\frac{1}{p} - \frac{1}{2r}} + \frac{\varepsilon}{2}. 
\end{align*}
so once again we can carry out the same argument, this time with the help of (II-ii) from \cite{Miyachi}$, \widehat{f}(\xi)= (1-\vartheta_0(\xi)) |\xi|^{-\lambda_1}$ and $\widehat{g}(\xi)= (1-\vartheta_0(\xi)) |\xi|^{-\lambda_2}$. We conclude that $f \in L^p$ and $g \in L^q$ but $B(f,g) \not\in L^r$.

Finally turning to the case of bilinear operators of the form
\begin{equation*}
T(f,g)(x)= \int_{\mathbb{R}^{2n} } a(\xi,\eta)\widehat{f}(\xi) \, \widehat{g}(\eta)\, e^{ix\cdot(\xi+\eta)}\, e^{ i \phase(\xi) - i \phase(\eta)+i \varphi_0 (\xi+\eta)} \, \dd\xi \, \dd\eta,
\end{equation*}
we observe that
$T(f,g)(x)= e^{i\varphi_0(D)} (B(f,g))(x),$ with $B(f,g)$ as above. Therefore the unitarity of the operator $e^{i\varphi_0(D)}$ on $L^2$ yields that the boundedness of $T$ from $L^p \times L^q \to L^2$ is equivalent to the $L^p \times L^q \to L^2$ boundedness of $B(\cdot, \cdot)$, and the discussion above establishes the sharpness of the parameters involved, in the case $1/p + 1/q=1/2$ and $r=2$.

\section{Basic Vector-Valued and Carleson Estimates for Oscillatory Integral Operators} \label{section main linear estim}
Before proceeding to the boundedness results, we need the following lemma, which was proved in the case of FIOs in \cite{RRS2}. We also include the proof, both for the sake of completeness and for later reference.
 
\begin{lemma}\label{lem:Wase-Lin}
Let $\vth \colon \R^n \to \R$ be a positive, radial, radially decreasing, smooth cut-off function which satisfies $\vth(\xi) = 1$ if $|\xi| \leq 1$ and $\vth(\xi) = 0$ if $|\xi| \geq 2$ $($as defined in \emph{Definition~\ref{def:LP}}$)$, and set $\th_k(\xi) := \vth(2^{3-k}\xi)$. Furthermore let $\omega_k(\xi)$ be a bump function equal to one on the support of $\th_k$. Now assume that \[
    s>0, \qquad s_c= \min(s,1),\qquad n/(n+s_c)<p\leq \infty,\qquad  m=-ns\abs{\frac{1}{p}-\frac{1}{2}},
\]
and for a fixed but arbitrary vector $u\in \R^n$ set
\[
    b(k,\xi) := 2^{km}\omegaF_k(\xi)\quad\mbox{and} \quad \widehat{P_k^{u}(g)}(\xi):=\th_k(\xi)e^{i2^{-k} \xi\cdot u}\widehat{g}(\xi).
    \]

If $\varphi$ is a phase function of order $s$, then one has 
\begin{equation}\label{eq:hp_composition}
\sup_k \norm{(P_{k}^{u} \circ T^{\phase}_{b})(f)}_{h^p}\lesssim \norm{f}_{h^p}, 
\end{equation}
and for $n\geq 1$ one also has for $m=-ns/2$
\begin{equation}\label{eq:uniform_bounddness}
    \sup_k \norm{(P_{k}^{u} \circ T^{\phase}_{b})(f)}_{L^\infty}\lesssim \norm{f}_{\bmo} \quad \mbox{and} \quad \sup_k \norm{(P_{k}^{u} \circ T^{\phase}_{b})(f)}_{h^1}\lesssim \norm{f}_{L^1}.
\end{equation}
The same conclusion holds for \emph{FIOs}, that is, when $s=1$ and $\phase$ is positively homogeneous of degree one. In that case \eqref{eq:hp_composition} is valid for  $m=-(n-1)\abs{\frac{1}{p}-\frac{1}{2}}$ and $n/(n+1)<p\leq \infty,$ and \eqref{eq:uniform_bounddness} is valid when $n\geq 2$ and $m=-(n-1)/2$.
\end{lemma}
\begin{proof}
The proof of \eqref{eq:hp_composition} follows from the fact that the amplitude of $P_{k}^{u} \circ T^{\phase}_{b}$ is in $S^m$ uniformly in $k$. 

In order to establish the first inequality in \eqref{eq:uniform_bounddness}, we write $b=b^\flat+b^\sharp$ where
\begin{equation}\label{decomposition of bj}
    b^\flat(k,\xi)=b(k,\xi)(1-\lambda(\xi)),\quad \mbox{and}\quad b^\sharp(k,\xi)=b(k,\xi)\lambda(\xi).
\end{equation}
and $\lambda$ is a smooth function that vanishes in a neighbourhood of the origin and equal to one outside a larger neighbourhood of the origin. Now since $m\leq 0$ and $1-\lambda$ is a low frequency cut-off, one can essentially throw away the $\omegaF$ in the definition of $b$ which would then make $b^\flat$ equal to $2^{km}$. Then by the kernel estimates for the OIOs with amplitude $b^\flat$ (see e.g. Lemma \ref{low freq lemma llf 1}), for $f\in \bmo$ we have that 
\[
        \norm{P_k^{u} T_{b^\flat}^{\phase}(f)}_{L^\infty}\lesssim \norm{T_{b^\flat}^{\phase}(f)}_{L^\infty}\lesssim \norm{(1-\lambda)(D)f}_{L^\infty}\lesssim \norm{f}_{\bmo}.
\]

In order to ameliorate $(P_{k}^{u} \circ T^{\phase}_{b^\sharp})(f)$ so that we can better understand its action on $\bmo$ functions, we employ an argument from \cite{Monster}*{page 27}.
According to that argument, for $n\geq 1$ and $m=\frac{-ns}{2}$, one introduces
the operator
$R_{k} (G)(x)=\int K_{k}(x-y) G(y)\dd y$, with
\begin{equation*}
	K_{k}(z)=\sum_{\kappa\leq j\leq k}2^{j m} \Psi\brkt{2^{k-j}z} 2^{n(k-j)},\quad \textrm{for some}\,\,\, \kappa
\end{equation*}
and 
\begin{equation}\label{t_eq:new_psi_2}
	\widehat{\Psi}(\eta)=\widehat{\psi}(\eta)^2\abs{\eta}^{m}=:\widehat{\psi}(\eta) \widehat{\tilde{\Psi}}(\eta),
\end{equation}
where $\widehat{\psi}$ is smooth, radial and positive with
\[
	\supp \widehat{\psi}\subset \set{\xi:\, 2^{-2}\leq \abs{\xi}\leq 1},
\]
and
\[
	\sum_{j\in \Z} \widehat{\psi}(2^{j} \eta)^2=1\qquad {\text{for any}\, \eta\neq 0}.
\]
Moreover, by \cite[Lemma~4.8]{Monster}, the kernel $K_k$ has the following properties:
\begin{equation*}
	\int K_k(z) \dd z=0;
\end{equation*}
and for each $0<\delta<\frac{ns}{2}$ the estimates
\begin{equation*}
	\abs{K_k(x-y)}\lesssim {2^{kn}}{\brkt{1+\frac{\abs{x-y}}{2^{-k}}}^{-n-\delta}}
\end{equation*}
and
\begin{equation*}
	\abs{K_k(x-y)-K_k(x-y')}\lesssim {2^{k(n+1)}}\abs{y-y'}
\end{equation*}
hold for all $x,y,y'\in \R^n$ and $k \in \Z$. Therefore the operator $R_k$ satisfies
\begin{equation*}
	\sup_{k\in\Z}\norm{R_k f}_{L^q}\lesssim \norm{f}_{L^q},\qquad 1\leq q< \infty,
\end{equation*}
and
\begin{equation*}\label{t_eq:BMO_bound}
\sup_{k\in\Z}\norm{R_k f}_{L^\infty}\lesssim \norm{f}_{\mathrm{BMO}}.
\end{equation*}
The consequence of the above discussion is that we can write
\begin{equation}\label{defn: Rk}
  R_k= \sum_{\kappa \leq j\leq k} Q_j 2^{(k-j)m}
\end{equation}
and $Q_j (D):= \widehat{\Psi}(2^{-j} D)$,  which enables one to replace $(P_{k}^{u} \circ T^{\phase}_{b^\sharp})(f)$ by $P_k^{u} \circ R_k \circ T^{\phase}_{\gamma}(f)$, for $n\geq 1$, where $\gamma(\xi):= \lambda(\xi)|\xi|^m \in S^{-ns/2}$.

Using the $\mathrm{BMO}$--$L^\infty$ boundedness above, the global $\bmo$-boundedness of OIOs with amplitudes in $S^{-ns/2}$ (i.e.\ Theorem \ref{linearhpthmoio}) and the $L^\infty$-boundedness of $P_k^u$, all together yields that 
\[
\sup_{k}\norm{P_k^{u} T_{b^\sharp}^{\phase}(f)}_{L^\infty} = \sup_{k}\norm{P_k^{u} \circ R_k \circ T^{\phase}_{\gamma}(f)}_{L^\infty}\lesssim \norm{\lambda(D)f}_{\mathrm{BMO}}\leq \norm{f}_{\bmo}.\qedhere
\]
\end{proof}

Another useful tool in our analysis is the following lemma.
\begin{lem}\label{Lu-Be} 
Let 
\[
    s>0, \quad n\geq 1,\quad n/(n+ s_c)<p< \infty, \quad p\neq 2 \quad \mbox{and} \quad m=-ns\abs{\frac{1}{p}-\frac{1}{2}}.
\]
{Assume that $ b(k,\xi)$ and $P_k^{u}$ are given by the same expressions as in \emph{Lemma~\ref{lem:Wase-Lin}}}.
Then for an \emph{OIO} $T^{\varphi}_b$ one has
\[
    \Big\Vert\Big(
\sum_{k=0}^{\infty}\, |P_k^{u} T_{b}^{\phase}(f)|^2\Big)^{1/2}\Big\Vert_{L^p} \lesssim \norm{f}_{h^p}.
\]
For an \emph{FIO} $T^{\varphi}_b$ the same result is valid under the conditions that  $n\geq 2,$ $m=-(n-1)\abs{\frac{1}{p}-\frac{1}{2}}$ and $n/(n+1)<p\leq \infty$ {with $p\neq 2$}.
\end{lem}
\begin{proof}
We only prove the result for the case of OIOs since the corresponding proof for FIOs is carried out in a similar manner. Observe that $P_k^{u} T_{b}^{\phase}$ is an oscillatory integral with amplitude $$e^{i2^{-k} \eta\cdot
u}2^{km}\widehat{\omega}_k(\eta)\widehat{\theta}(2^{-k} \eta)$$ and phase function $x\cdot\eta+ \varphi (\eta).$ 
One can also write the amplitude as
\begin{equation*}\label{amplitude of T2k}
\begin{split}
e^{i2^{-k} \eta\cdot
u}2^{km}&\widehat{\omega}_k(\eta)\widehat{\theta}(2^{-k} \eta)\\
&= \lambda(\eta) e^{i2^{-k} \eta\cdot
u}2^{km}\widehat{\omega}_k(\eta)\widehat{\theta}(2^{-k} \eta)+(1-\lambda(\eta)) e^{i2^{-k} \eta\cdot
u}2^{km}\widehat{\omega}_k(\eta)\widehat{\theta}(2^{-k} \eta)\\
&:= \alpha_k +\beta_k,  
\end{split}
 \end{equation*}
where $\lambda$ is the high frequency localisation introduced in \eqref{decomposition of bj}.

We first consider the case of $\Big\Vert\Big(
\sum_{k=0}^{\infty}\, |P_k^{u}T^{\varphi}_{\alpha_k} (f)|^2\Big)^{1/2}\Big\Vert_{L^p} .$ Replacing $T^{\varphi}_{\alpha_k}$ with $P_{k}^{u} \circ R_k\circ T^{\phase}_{\gamma}$, with $\gamma\in S^{m}_{1,0}$,  matters reduce to proving the desired boundedness for
\[
   ( \sum_{k\geq0} |(P_{k}^{u} \circ R_k\circ T^{\phase}_{\gamma})f|^{2})^{1/2},
\]
where $T_\gamma^{\varphi}$ and $R_{k}$ are as in Lemma \ref{lem:Wase-Lin}. 
Now if we introduce a smooth cut-off function ${\chi}$ such that $R_{k}= R_{k} (1-{\chi}(D))$ and using Theorem \ref{linearhpthmoio} for OIOs (or Theorem \ref{linearhpthmfio} in the case of FIOs), it is enough to prove 
\begin{equation}\label{eq:Rk quasibanach}
\brkt{\int {\brkt{\sum_{k\geq 0} \abs{P_{k}^{u} R_{k} G(x)}^2}^{p/2}}\dd x}^{1/p} \lesssim \Vert G\Vert_{h^p}.
\end{equation}
At this point, for the sake of simplicity of the notation, we replace $P_{k}^{u}$ by $P_{k}$ in what follows. This modification will not cause any problems since the difference between the two operators only lies in a harmless factor $e^{i({\cdot})\cdot u}$. Observe now that using the integral representation of $R_k$, one has
\[
	P_{k} R_{k} G(x)=\sum_{\kappa\leq j\leq k} 2^{jm} \int \big(\theta_{k}*\tilde{\Psi}_{2^{j-k}}\big)(y) \big(\psi_{ 2^{j-k}}*G\big)(x-y)\dd y,
\]
where $\tilde{\Psi}$ is defined in \eqref{t_eq:new_psi_2} and $\tilde{\Psi}_{(\cdot)}(x):= (\cdot)^{-n} \tilde{\Psi}(\frac{x}{(\cdot)}),$ and $\psi_{ 2^{j-k}}$ is defined in a similar way.
 Then for any $\nu>0$ (to be later determined)
\[
	\abs{P_{k} R_{k} G(x)}\leq \sum_{\kappa\leq j\leq k} 2^{jm} \left( \int \abs{\theta_{k}*\tilde{\Psi}_{2^{j-k}}(y)}\brkt{1+\frac{\abs{y}}{2^{j-k}}}^\nu\dd y\right) \, \mathcal{M}_{\nu, 2^{k-j}}(\psi_{2^{j-k}}\ast G)(x),
\]
where $\mathcal{M}_{\nu, 2^{k-j}}$ is the Peetre maximal function as defined in \eqref{Peetremax}.

Now by \eqref{hl bounds peetre} we have
\[
	\mathcal{M}_{\nu, 2^{k-j}}(\psi_{2^{j-k}}\ast G)(x)\lesssim \mathcal{M}_{n/\nu}(\psi_{2^{j-k}}*G),
\]
for any $x\in \R^n$. Moreover by fairly standard estimates for convolution-type integrals one has for any $N>n+\nu$
\[
	 \abs{\theta_{k}*\tilde\Psi_{2^{j-k}}(y)}\lesssim  ({2^{-k}}\max\brkt{2^j,1})^{-n}
	 \brkt{1+\frac{\abs{y}}{2^{-k}\max\brkt{2^j,1}}}^{-N},
\]
which in turn implies that
\[
	 \sup_{\kappa\leq j\leq k}\int \abs{\theta_{k}*\tilde\Psi_{2^{j-k} }(y)}\brkt{1+\frac{\abs{y}}{2^{j-k}}}^\nu\dd y\lesssim 1.
\]
Thus, we have the pointwise inequality
\begin{equation*}\label{eq:Pointwise_estiamte_RK}
	\abs{P_{k} R_{k} G(x)}\lesssim \sum_{\kappa\leq j\leq k}2^{jm}\left[\mathcal{M}\brkt{\abs{\psi_{2^{j-k}}*G}^{\frac{n}{\nu}}}(x)\right]^\frac{\nu}{n}.
\end{equation*}

Therefore, for any $q>\max\brkt{\frac{n}{\nu},1}$
\begin{equation*}
\begin{split}
	\brkt{\sum_{k\geq 0} \abs{P_{k} R_{k} G(x)}^{q}}^{1/q}&\lesssim \sum_{j=\kappa}^\infty 2^{jm}\brkt{ \sum_{k\geq j}\left[\mathcal{M}\brkt{\abs{\psi_{2^{j-k}}*G}^{\frac{n}{\nu}}}(x)\right]^\frac{q\nu}{n}}^{\frac{1}{q}}\\
	&\leq C_{m} \brkt{ \sum_{k\geq 0} \left[\mathcal{M}\brkt{\abs{\psi_{{k}}*G}^{\frac{n}{\nu}}}(x)\right]^\frac{q\nu}{n}}^{\frac{1}{q}} \\&= C_{m} \brkt{ \sum_{k\geq 0} \left[\mathcal{M}_{n/\nu}\brkt{\psi_{{k}}*G}(x)\right]^q}^{\frac{1}{q}}
	\end{split}
\end{equation*}
where $C_{m}=\sum_{j=\kappa}^\infty 2^{jm}<+\infty$.

Hence, for any $p>\frac{n}{\nu}$ the Fefferman-Stein's  estimate \eqref{FSvector}
yields that 

\[
	\norm{\brkt{\sum_{k\geq 0} \abs{P_{k} R_{k} G}^{q}}^{1/q}}_{L^p(\R^n)}\lesssim \left[\int \brkt{\sum_{k\geq 0} {\abs{\psi_{{k}}*G(x)}}^q}^{\frac{p}{q}}\dd x\right]^{\frac{1}{p}}.
\]
Finally the last term is equal to
\[
	\left[\int \brkt{\sum_{k\geq 0} {\abs{\psi_{{k}}*G(x)}}^q}^{\frac{p}{q}}\dd x\right]^{\frac{1}{p}}	\lesssim \norm{G}_{F^0_{p,q}}.
\]
 Taking $q=2$ and $\nu > n/p$, and using Definition \ref{def:Triebel}, we obtain
\[
	\brkt{\int {\brkt{\sum_{k\geq 0} \abs{P_{k} R_{k} G(x)}^2}^{p/2}}\dd x}^{1/p}\lesssim \norm{G}_{h^{p}}.
\]
This proves \eqref{eq:Rk quasibanach}.\\

Now to treat $\Big\Vert\Big(\sum_{k=0}^{\infty}
\, |P_k^{u}T^{\varphi}_{\beta_k} (f)|^2\Big)^{1/2}\Big\Vert_{L^p},$  we observe that by an argument similar to the proof of Lemma \ref{lem:Wase-Lin} (that is to say, essentially use \eqref{eq:hp_composition}) we have
\begin{equation}
\begin{split}
 \Big\Vert\Big(\sum_{k=0}^{\infty}
\, |P_{k}^{u}T^{\varphi}_{\beta_k} (f)|^2\Big)^{1/2}\Big\Vert_{L^p}^{p}
&\lesssim  \Big\Vert\Big(\sum_{k=0}^{\infty}\, 2^{km} |P_{k}^{u} T_{1-\lambda}^{\varphi} f|^2\Big)^{1/2}\Big\Vert_{L^p}^p\\
&\lesssim 
\sum_{k=0}^{\infty}\, 2^{pkm/2}     \norm{P_{k}^{u} T_{1-\lambda}^{\varphi_2} f}_{L^p}^p \Big)\\
& \lesssim  \brkt{\sum_{k=0}^{\infty}  2^{pkm/2}} 
\sup_{k\geq 0} \norm{P_{k}^{u} T_{1-\lambda}^{\varphi_2} f}_{L^p}^p\\
& \lesssim  \norm{f}_{h^p}^{p}.\qedhere
\end{split}
\end{equation}
\end{proof}
\begin{rem}\label{PDErem}
A re-examination of the proofs of \emph{Lemmas} \emph{\ref{lem:Wase-Lin}} and \emph{\ref{Lu-Be}} reveals that if $b(k, \xi)= 2^{km_0}\omega_{k}$ with $m_0<0$ and $1<p<\infty,$ $m(p)=-ns\Big|\frac{1}{p}-\frac{1}{2}\Big|$ then one has
\[
\Big\Vert\Big(
\sum_{k=0}^{\infty}\, |P_k^{u} T_{b}^{\phase}(f)|^2\Big)^{1/2}\Big\Vert_{L^p} \lesssim \norm{f}_{H^{m_0-m(p),p}}.
\]
where $H^{s,p}= F^{s}_{p,2}$ is the $L^p$-based Sobolev space.
\end{rem}
For the boundedness of  multilinear OIOs with target spaces $L^2$ or $\bmo$, we need the following result about oscillatory integrals giving rise to Carleson measures, whose counterpart in the case of FIOs was proven in \cite{RRS2}. The proposition below doesn't require any homogeneity from the phase function as in the case of FIOs.
\begin{prop}\label{lem:smallcarlnorm}
Let $s>0$, $d \in S^{-ns/2}$, $u\in \R^n$ and let 
\[
    Q^u_{k} f(x)= \frac{1}{(2\pi)^n} \int \psiF_k(\xi)e^{i2^{-k} \xi\cdot u}\widehat{f}(\xi) e^{ix\cdot \xi} \dd \xi,
\]
where  $k ≥ k_0 ∈ ℤ$ and
\begin{equation}\label{defn of psik}
    \psiF_k(\xi)^2 := \vth(2^{-1-k}\xi)^2 - \vth(2^{2-k}\xi)^2 ,
\end{equation}
and $\vartheta$ is as in \emph{Lemma \ref{lem:Wase-Lin}}. Then if $\varphi$ is a phase function of order $s>0$ and $f \in \bmo$ one has that
\m{
    \dd \mu_k(x,t) = \sum_{ℓ=0}^\infty |(Q_{k+ℓ}^{u} \circ T^{\phase}_{d})(f)(x)|^2 \delta_{2^{-ℓ}}(t) \dd x
}
is a Carleson measure with Carleson norm bounded by $C_ε 2^{-εk}\|f\|_\bmo^2$. Here, for any $δ ∈(0,1)$, $\varepsilon$ is given by $\min(ns/2, nδ)$.  
\end{prop}

\begin{proof}
Since we can write $Q^u_{k+ℓ}∘T_d^φ = Q^u_{k+ℓ}∘T_d^φ∘\tilde Q^u_{k+ℓ}$, where $\tilde Q^u_{k+ℓ} : \bmo → L^∞$ uniformly in $k$, we first consider the case of $f∈L^∞$. Also for simplicity of the exposition we set $u=0$ in what follows.

Now since the operator $Q_{k+ℓ} ∘ T^φ_d$ is essentially the $(k+\ell)$-th component of the Littlewood--Paley decomposition of the operator $T^φ_d$, setting $j = k + ℓ ≥ k_0$  we carry out a second microlocalisation of $Q_{j} ∘ T^φ_d$ in the following way.

Take a non-negative real number $μ$, to be fixed later, and for each $j$, fix $O(2^{nμj})$ vectors $ξ_j^ν$, $ν = 1, \ldots, O(2^{nμj})$, distributed evenly in $\supp \psi_j$. Let $\{ρ_j^ν\}_ν$ be a family of smooth functions, where $\supp ρ_j^ν$ is a ball of radius $2^{(1-μ)j}$ centred at $ξ_j^ν$, chosen in such a way that the supports of $\{ρ_j^ν\}_\nu$ cover $\supp \psi_j$. One may for example take a smooth bump function $β$ supported in a ball of radius $1$ about the origin and from this form $ρ_j^ν(ξ) = β(2^{(μ-1)j}(ξ-ξ_j^ν))/∑_κβ(2^{(μ-1)j}(ξ-ξ_j^κ))$.

It is clear that these cut-off-functions satisfy
\[
    |∂^αρ_j^ν(ξ)| ≤ C_α 2^{|α|(μ-1)j}.
\]
With this partition of unity, we may therefore write the integral kernel of $Q_j∘T_d^φ$ as $K_j(x,y) = ∑_ν K_j^ν(x,y)$, with
\[
    K_j^ν(x,y) = ∫ d (ξ) ρ_j^ν(ξ)ψ_j(ξ)e^{i(x-y)·ξ + iφ(ξ)}\ddd ξ.
\]
In order to get desired estimates for the kernel, we rewrite the phase of this integral as
\begin{align*}
    (x-y)·ξ + φ(ξ) &= (x - y + ∇φ(ξ_j^ν))·ξ + h_j^ν(ξ),\\[.5em]
    \text{with}\qquad h_j^ν(ξ) &= φ(ξ) - ∇φ(ξ_j^ν)·ξ,
\end{align*}
which in turn yields
\[
    K_j^ν(x,y) = ∫b_j^ν(ξ) e^{i(x - y + ∇φ(ξ_j^ν))·ξ}\ddd ξ, 
\]
where $b_j^ν(ξ) = d (ξ)ρ_j^ν(ξ)ψ_j(ξ)e^{ih_j^ν(ξ)}$. The mean-value theorem then yields that 
$∂_ih_j^ν(ξ) =  ∇∂_iφ(η)·(ξ-ξ_j^ν)$ for some $η$ on the line segment between $\xi$ and $\xi^{\nu}_j$. On $\supp ψ_j\,ρ_j^ν$, we therefore have from (\ref{eq:simplified phase}) that
\begin{align*}
    |∂^α h_j^ν(ξ)| &\lesssim\begin{cases} 2^{(s-μ-1)j} &|α|=1\\ 2^{(s-|α|)j}& |α|>1.
    \end{cases}
\end{align*}
If we take $μ ≤ s/2$, the worst terms of $|∂_ξ^αe^{ih_j^ν(ξ)}|$ are hence bounded by a constant times $2^{(s - μ -1)|α|j}$.

With these estimates at hand, we find that on the support of $ψ_j\,ρ_j^ν$, 
\begin{align*}
        |∂^α b_j^ν(ξ)| &≤ C_α ∑_{∑α_ℓ=α}\big|∂^{α_1}d(ξ)∂^{α_2}ρ_j^ν(ξ)∂^{α_3}ψ_j(ξ)∂_ξ^{α_4}(e^{ih_j^ν(ξ)})\big|\\
                &≤ C_α ∑_{∑α_ℓ = α}2^{(-ns/2 - |α_1 + α_3| + (μ-1)|α_2| + (s-μ-1)|α_4|)j} ≤ C_α 2^{(-ns/2 + (s/2-1)|α|)j},
\end{align*}
where we have fixed $μ$ to be the optimal $μ = s/2$. For later convenience, we define $λ = s/2 - 1$.

We are now ready to take on the Carleson norm estimates. To that end, we fix a ball $B$ of radius $r<1$ and centre $x_0$.
Let then $τ ∈ (0,1)$ be given by
\[
    τ = \begin{cases}1 - s/2 & \text{if }s < 2\\ 1 - δ & \text{otherwise,}\end{cases}
\]
where $δ ∈ (0,1)$ is arbitrary, and
let $R_ν$ be the ball of radius $ 2· 2^{(λ+τ)(j-k_0)}r^τ$ and centre $x_0 + ∇φ(ξ_j^ν)$. 
Clearly then,
\[
    Q_jT_d^φf(x) = ∑_ν S^\nu_j(χ_{R_ν}f)(x) + ∑_ν∫_{R_ν^c} K_j^ν(x,y)f(y)\dd y,
\]
where $S^\nu_j$ is the operator with kernel $K^\nu_j$.

For the parts inside the balls $R_ν$, we use that $|ψ(2^{-j}ξ)d(ξ)| \lesssim 2^{-nsj/2}$, and hence $S_j^ν$ is bounded $L^2 → L^2$ with operator norm estimated by $2^{-nsj/2}$.
Using this and the fact that the symbols have almost disjoint support -- that is, with a finite number of overlaps -- we find that for each $ν$,
\begin{align*}
     ∫_B \Big|\sum_\nu S^ν_j(χ_{R_ν} f)(x) \Big|^2 \dd x &\lesssim 
     ∫\Big|\sum_\nu S^ν_j(χ_{R_ν} f)(x) \Big|^2 \dd x\lesssim \sum_\nu 2^{-nsj}\|χ_{R_ν} f\|^2_{L^2}\\ &≤  
      \sum_\nu 2^{-nsj} |R_ν| \|f\|^2_{L^∞}
      \lesssim (2^jr)^{(τ-1)n} |B| \|f\|^2_{L^∞}.
\end{align*}

To find a similar estimate for the parts outside $R_ν$ we start by noting that the triangle inequality and the fact that $λ + τ ≥ 0$ and $j≥k_0$ yield that for $r\leq 1$, any $x∈B$ and any $y$ with $y + x_0 + ∇φ(ξ^ν_j) ∈ R_ν^c$ we have
\[
    |x - x_0 - y| ≥ |y| - r ≥  \frac12|y| + 2^{(λ + τ)(j-k_0)j}r^τ + r ≥ \frac12|y| \gtrsim 2^{(λ +τ)j}r^τ.
\]
We therefore have for any $x∈B$ and non-negative integer $N$ that
\begin{align*}
    ∫_{R_ν^c}|K_j^ν(x,y)|\dd y &= ∫_{R_ν^c} \Big|∫ b_j^ν(ξ)\,\Big(\frac{(x+∇φ(ξ_j^ν)-y)·∇_ξ}{|x + ∇φ(ξ_j^ν)- y|^2}\Big)^N e^{i(x + ∇φ(ξ_j^ν) - y) ·ξ} \ddd ξ \Big|\dd y\\
        &= ∫_{R_ν^c} \Big|∫ e^{i(x + ∇φ(ξ_j^ν) - y) ·ξ} \Big(\frac{(x+∇φ(ξ_j^ν)-y)·∇_ξ}{|x + ∇φ(ξ_j^ν)- y|^2}\Big)^Nb_j^ν(ξ) \ddd ξ \Big|\dd y\\
        &\lesssim ∫_{R_ν^c}\frac{2^{-nsj/2 + Nλj}|\supp ρ_j^ν|}{|x + ∇φ(ξ_j^ν)- y|^N}\dd y\\
        &= ∫_{y + x_0 + ∇φ(ξ_j^ν) \in R_ν^c}\frac{2^{-nsj/2 + Nλj}|\supp ρ_j^ν|}{|x -x_0 - y|^N}\dd y\\
        &\lesssim ∫_{|y| ≥ 2^{(λ + τ)j}r^τ}\frac{2^{-nsj/2 + (N-n)λj}}{|y|^N}\dd y 
        \ \lesssim\ 2^{-nsj/2}(2^jr)^{(n-N)τ}.
\end{align*}
Now choose $N$ large enough to make $2(N - n)τ ≥ n(1 - τ) = : ε$. Note that from the definition of $τ$, we have that $ε = \min(ns/2, nδ)$, where $δ ∈ (0,1)$ is arbitrary.
Combining with the result from the part from inside  $R_ν$ and summing over the $O(2^{nsj/2})$ balls, this then yields that $∫_B |Q_jT_d^φ(f)(x)|^2\dd x \lesssim (2^jr)^{-ε}|B|\|f\|_{L^∞}^2$. Hence
\begin{align}\label{small_carleson}
    ∫_{B×[0,r]}|\mathrm dμ_k(x,t)| &= ∑_{2^{-ℓ}≤r}∫_B |Q_{k+ℓ}∘T_d^φ∘\tilde Q_{k+ℓ}(f)(x)|^2\dd x\\
\nonumber   &\lesssim∑_{2^{-ℓ}≤r}(2^{k+ℓ}r)^{-ε}|B|\|\tilde Q_{k+ℓ}(f)\|_{L^∞}^2 \\
\nonumber   &\lesssim 2^{-εk}∑_{2^{-ℓ}≤r}(2^{-ε})^ℓr^{-ε}|B|\|f\|_\bmo^2 \lesssim 2^{-εk}|B|\|f\|_\bmo^2,
\end{align}
which shows the requested Carleson estimate for balls of radius smaller than $1$.

Now if the radius $r$ of $B$  is larger than one, then we cover $B$ by balls $B_j$ of radius $1/2$, observing that there are $O(r^n)$ such balls needed for this covering. Furthermore we observe that for $r>1$, \eqref{small_carleson} yields that
\ma{
    \int_{B\times [0,r]} |\dd\mu_k(x,t)| &=\int_{B\times [0,1]} |\dd\mu_k(x,t)|\leq 
    \sum_{O(r^n)} \int_{B_j\times [0,1]} |\dd\mu_k(x,t)|\\
&\lesssim  \sum_{O(r^n)} 2^{-\varepsilon k} 2^{-n}\norm{f}_{\bmo}^2\lesssim 2^{-\varepsilon k} |B| \norm{f}_{\bmo}^2. \qedhere
}
\end{proof}

\section{Frequency decomposition of the oscillatory integral operator}
\label{sec:freqdecom}

Following the method in \cite{RRS2} for the decomposition of the amplitude $\sigma(x, \Xi)\in S^m(n,N)$, we reduce the problem of regularity of $T^\Phase_\sigma$ into considering three frequency regimes: When $\Xi$ lies inside a compact set; when one component of $\Xi= (\xi_1, \dots, \xi_N)$ dominates the others; and when two fixed components of $(\xi_1, \dots, \xi_N)$ are comparable to each other. In what follows we only describe the aspects of the amplitude decomposition which are crucial to the later sections of the paper. For the remaining details, we refer the reader to \cite{RRS2}.

Here and in all that follows we take $N>1$. First we define the component of $\sigma$ with frequency support contained in a compact set. We introduce a cut-off function $\chi \colon \R^{nN} \to \R$, such that $\chi(\Xi) = 1$ for $|\Xi| \leq 1/8$ and $\chi(\Xi) = 0$ for $|\Xi| \geq 1/4$ and define
\begin{equation}\label{defn of sigma0}
    \sigma_0(x,\Xi) = \chi(\Xi)\,\sigma(x,\Xi).
\end{equation}
To define the components of $\sigma$ where one frequency dominates all the others, we construct a cut-off function $\nu \colon \R^{nN} \to \R$ such that $\nu(\Xi) = 0$ for $|\xi_1| \leq 32\sqrt{N-1}\,|\Xi'|$ and $\nu(\Xi) = 1$ for $64\sqrt{N-1}\,|\Xi'| \leq |\xi_1|$, where $\Xi' := (\xi_2,\ldots,\xi_N)$. 
This can be done by taking $\Lambda\in \mathcal{C}^\infty(\R)$ such that $\Lambda(t)=1$ if $t\leq c_1$ and $\Lambda(t)=0,$ if $t\geq c_2$ for two suitably chosen real numbers $0<c_1<c_2<1$.

Define 
\begin{equation}\label{biglambda}
    \nu(\Xi)=1-\Lambda\brkt{\frac{\abs{\xi_1}^2}{\abs{\Xi}^2}}\in\mathcal{C}^\infty(\R^{nN}\setminus 0).
\end{equation}
Now given $j=1,\ldots N$ we define $\Xi'_j:=(\xi_1,\ldots, \xi_{j-1},\xi_{j+1},\ldots ,\xi_N)$ and 
\[
    \nu_j(\Xi):=\nu(\xi_j,\Xi'_j),
\]
for all $\Xi\in \R^{nN}$. We then define the component of $\sigma$ for which $\xi_j$ dominates the other frequency components to be
\begin{equation}\label{defn:sigmaj}
    \sigma_j(x,\Xi) = (1-\chi(\Xi))\,\nu_j(\Xi)\,\sigma(x,\Xi), \quad \mbox{for $j=1,\ldots N$.}
\end{equation}

What remains of $\sigma$ will be split into functions on whose support two frequency components are comparable (see \cite{RRS2} pages 22--23 for the details). Thus $\sigma$ can be finally decomposed as
\begin{equation*}\label{main amplitude decomp}
    \sigma(x,\Xi) = \sigma_0(x,\Xi) + \sum_{j=1}^N \sigma_j(x,\Xi) + \sum_{j\neq k} \sigma_{j,k}(x,\Xi),
\end{equation*}
where $\sigma_0$ has compact $\Xi$-support, $\abs{\xi_j}$ dominates $\abs{\Xi}$ on the $\Xi$-support of $\sigma_j$, and $\abs{\xi_j} \approx \abs{\xi_k}$ on the $\Xi$-support of $\sigma_{j,k}$. More specifically, $\sigma_{j,k}$ and $\sigma_j$ are supported away from the origin, and
\begin{equation}\label{the key on the support of sigma1}
   c\abs{\xi_j}^2\geq \abs{\Xi}^2  
\end{equation}
on the $\Xi$-support of $\sigma_j$, for a suitably chosen $c>1$.

One can also check that if $\sigma \in S^m(n,N)$ then $\sigma_j$ and $\sigma_{j,k}$ are also in $S^m(n,N)$ for all $j,k=1,\dots,N$ and $\sigma_0 \in S^\mu(n,N)$ for all $\mu \in \R$.\\

We shall now proceed by giving  explicit representations for the multilinear OIOs  $T_{\sigma_{0}}^\Phase$, $T_{\sigma_{1}}^\Phase$ and $T_{\sigma_{1,2}}^\Phase$, which as will be clarified in Section \ref{endpoint cases}, are the prototypes of the operators for which the boundedness results will be established here. Moreover, the boundedness of $T_{\sigma}^\Phase$ can be reduced to the boundedness of these three types of operators. However further reductions are needed to make the representations of the aforementioned operators amenable to the vector-valued- and maximal-function-based proofs that are utilised in this paper.


\subsection{Representation of $T_{\sigma_{0}}^\Phase$}
We note that by \eqref{defn of sigma0}, the support of $\sigma_0$ is in a fixed compact set. Therefore as was demonstrated in \cite[page 44]{RRS2} the operator $T^\Phase_{\sigma_{0}}$ can be written as 
\begin{equation}\label{final form of sigma0}
    T^\Phase_{\sigma_{0}} (f_1,\dots,f_N)(x) = \sum_{K\in\Z^{nN}} a_K(x)T_{\theta(\cdot/\sqrt{N})}^{\phase_0}\left(\prod_{j=1}^N T_{\theta}^{\phase_j}\circ\tau_{\frac{2\pi k_j}{L}}(f_j)\right)(x),
\end{equation}
where $\tau_h f(x) := f(x-h)$, $\theta \in \mathcal{C}_c^{\infty}(\R^n)$ and $a_K(x)$ is a smooth function satisfying \begin{equation}\label{derivative of a_K}
 |\partial^{\alpha}a_K(x)|\lesssim (1+\sum_{j=1}^{N} |k_j|^2)^{-M}   
\end{equation} for all $x\in \R^n$ and $M\geq 0$, with $K=(k_1,\dots,k_N)$.

\subsection{Representation of $T_{\sigma_{1}}^\Phase$}
Let $\vth$ be the function  introduced in Definition~\ref{def:LP} and recall or define

\begin{itemize}

\item $\th_k(\xi) := \vth(2^{3-k}\xi)$,

\item $\psiF_k(\xi)^2 := \vth(2^{-1-k}\xi)^2 - \vth(2^{2-k}\xi)^2$,

\item $\phiF_k(\xi)^2 := \vth(2^{-3-k}\xi)^2 - \vth(2^{4-k}\xi)^2$.

\end{itemize}

From the support properties of \(\sigma_{1}\), it follows that if \(\psi_{k}\left(\xi_{1}\right) \neq 0\) and \(\sigma_{1}(x, \Xi) \neq 0\)
then

$$
\left|2^{-k} \Xi_{1}^{\prime}\right| \leqslant \frac{\left|2^{-k} \xi_{1}\right|}{32 \sqrt{N-1}} \leqslant \frac{2^{-3}}{\sqrt{N-1}},
$$
which implies that \(\theta_{k}\left(\xi_{j}\right)=1\) for \(j=2, \ldots, N\), and one also has that
\nma{ineqforcut-off}{
\frac{1}{8} &\leq |2^{-k}(\xi_1+\dots+\xi_N)| < 8 \\
&\quad \mbox{which implies} \quad \phiF_k(\xi_1+\dots+\xi_N) = 1. 
}
Using these facts, there exists $k_0\in \Z$ (independent of $x$) such that we can write $T^\Phase_{\sigma_1}$ as
\begin{align*}\label{decomp12}
    &T^\Phi_{\sigma_1}(f_1,\ldots,f_N)(x) \\
    &= \int_{\R^{nN}} \sum_{k\geq k_0} \psiF_k(\xi_1)^2 \prod_{j=2}^N\th_k(\xi_j)^2 \phiF_k(\xi_1+\ldots+\xi_N)^2 \sigma_1(x,\Xi) \widehat{f_1}(\xi_1)\\&\times\prod_{j=2}^N\widehat{f_j}(\xi_j)\, e^{ix\cdot(\xi_1 +\dots+\xi_N)} e^{i\Phase(\Xi)} \ddd\Xi .
\end{align*}
See \cite[page 24]{RRS2} for the details of all these deductions.

We also introduce a high frequency cut-off $\chi_0$ that satisfies 
\begin{equation*}\label{eq:chi_0}
    \begin{cases}
     \chi_0(\xi) = 1,\qquad \mbox{for $|\xi| \geq 2^{k_0-4}$, and}\\
    \chi_0(\xi) = 0,\qquad \mbox{for $|\xi| \leq 2^{k_0-5}$,}
    \end{cases}
\end{equation*}
where $k_0$ can be chosen appropriately, and let $m_0,\dots,m_N$ be a (non-integer) partition of the decay $m$ of the amplitude $\sigma$, so that 
\begin{equation*}\label{defn of the total order}
   m = \sum_{j=0}^N m_j ,
\end{equation*} 
and $m_j = -ns|1/p_j - 1/2|$.
Based on these frequency cut-offs, we introduce the following localisation operators as well as amplitudes
\begin{align*}
\widehat{Q_k^0(f)}(\xi) &= \phiF_k(\xi)\widehat{f}(\xi), & b_0(\xi) &= |\xi|^{m_0}\chi_0(\xi), \\
\widehat{Q_k^{u_1}(f)}(\xi) &= |2^{-k}\xi|^{m-m_0-m_1}\psiF_k(\xi)e^{i2^{-k} \xi\cdot u_1}\widehat{f}(\xi), & b_1(\xi) &= |\xi|^{m_1}\chi_0 (\xi), \\
\widehat{P_k^{u_j}(f)}(\xi) &= \th_k(\xi)e^{i2^{-k} \xi\cdot u_j}\widehat{f}(\xi), & b_{j,k}(\xi) &= 2^{km_j}\omegaF_k(\xi),
\end{align*}
for $j=2,\dots,N$, $\omegaF_k(\xi) := \th_k(\xi/2)$ is the bump function introduced in Lemma~\ref{lem:Wase-Lin} equal to one on the support of $\th_k$.\\
{Also note that for any $m\leq 0$ the symbol
$2^{km}\omega_k(\xi)\in S^{m}$ uniformly in $k$, since
when $m\leq 0$ one has that $|2^{km}ω(2^{-k}ξ)| \lesssim 2^{km}\jap{2^{-k}ξ}^{m} ≤ \jap{ξ}^{m}$, since $ω$ is Schwartz, and moreover we also have that for any $N\geq 0$ and $|\alpha|>0$
$$|\partial^{\alpha}( 2^{km}ω(2^{-k}ξ))|\lesssim 2^{km} 2^{-k|\alpha|}(1+2^{-k} |\xi|)^{-N}\lesssim 2^{km} 2^{-k|\alpha|}2^{kN}(1+|\xi|)^{-N},$$
which by choosing $N= |\alpha|-m\geq 0,$ yields that 
$|\partial^{\alpha}( 2^{km}ω(2^{-k}ξ))|\lesssim \jap{\xi}^{m-|\alpha|}.$}

Using these operators one can show \cite[page 26]{RRS2} that for any $M\geq 0$, the operator $T_{\sigma_1}^{\Phi}$ can be written as
\begin{multline}\label{t_piecedecomp for sigma1}
  T_{\sigma_1}^{\Phi}(f_1,\ldots,f_N)(x) \\
 = \int \sum_{k\geq k_0}^\infty  M_{\mathfrak m}\circ T^{\phase_0}_{b_0} \circ P_{k}^0\left[(Q_{k}^{u_1} \circ T^{\phase_1}_{b_1})(f_1)\,\prod_{j=2}^N(P_{k}^{u_j} \circ T^{\phase_j}_{b_{j,k}})(f_j)\right](x)\\\times \frac{e^{i2^{-k}\Xi\cdot U}}{(1+|U|^2)^M} \, \ddd U,   
\end{multline}

where $M_{\mathfrak m}$ denotes the operator of multiplication by
{$\mathfrak m= \mathfrak m (k,x, U)$ with $U= (u_1, \dots, u_N)$, and $\mathfrak{m}$} is a smooth function depending on $\sigma_1$, with uniformly bounded derivatives of all orders.
It was shown in \cite[page 26]{RRS2} that boundedness of $T_{\sigma_{1}}^\Phase$ can been reduced to showing the boundedness of
\nma{the magic operator}{
    & {B}(f_1,\ldots,f_N)(x) \\
    &:={ \sum_{k\geq k_0 } \chi_0(2D)\,Q_{k}^0  \left[(Q_{k}^{u_1} \circ T^{\phase_1}_{b_1})(f_1)\prod_{j=2}^N(P_{k}^{u_j} \circ T^{\phase_j}_{b_{j,k}})(f_j)\right](x),}}
where the symbol of the high-frequency cut-off $\chi_0$ belongs to $S^0$.

\subsection{Representation of $T_{\sigma_{1,2}}^\Phase$}
With the same choice of $\psiF_k$, $\th_k$, $\chi_0$ and $\omega_k$ as above, and with a suitable choice of the integer $k_1$ and setting
\[
{\zeta}_{k}(\xi)^{2}:=\vartheta\big(2^{-k-k_{1}-2} \xi\big)^{2}-\vartheta\big(2^{3+k_{1}-k} \xi\big)^{2},
\]
it was demonstrated in \cite{RRS2} page 42, that for some $k_0\in \Z$ one has the representation
 
\[
\begin{split}
T_{\sigma_{1,2}}^{\Phi}\left(f_{1}, \ldots, f_{N}\right)(x) \\ =\int_{\mathbb{R}^{n N}} \sum_{k \geqslant k_{0}} \psi_{k}\left(\xi_{1}\right)^{2} {\zeta}_{k}\left(\xi_{2}\right)^{2} \sigma_{1,2}(x, \Xi) \chi_{0}\left(\xi_{1}\right) \widehat{f}_{1}\left(\xi_{1}\right) \times \\ \chi_{0}\left(\xi_{2}\right) \widehat{f}_{2}\left(\xi_{2}\right) \prod_{j=3}^{N} \theta_{k}\left(\xi_{j}\right)^{2} \widehat{f}_{j}\left(\xi_{j}\right) e^{i x \cdot\left(\xi_{1}+\cdots+\xi_{N}\right)+i \Phi(\Xi)} \mathrm{d} \Xi.
\end{split}
\]

Now we introduce the following localisation operators and amplitudes:
\begin{equation}
    \label{the d's}
\begin{aligned}
\widehat{P_k^0(f)}(\xi) &= \th_k(\xi)\widehat{f}(\xi), & d_{0}(\xi) &= 2^{km_0}\omegaF_k(\xi), \\
\widehat{Q_k^{u_1}(f)}(\xi) &= \abs{2^{-k}\xi}^{m-m_1-m_2}\psiF_k(\xi)e^{i2^{-k} \xi\cdot u_1}\widehat{f}(\xi), & d_1(\xi) &= |\xi|^{m_1}\chi_0 (\xi), \\
\widehat{Q_k^{u_2}(f)}(\xi) &= \psi_k(\xi) e^{i2^{-k} \xi\cdot u_2}\widehat{f}(\xi), & d_2(\xi) &= |\xi|^{m_2}\chi_0 (\xi), \\
\widehat{P_k^{u_j}(f)}(\xi) &= \th_k(\xi)e^{i2^{-k} \xi\cdot u_j}\widehat{f}(\xi), & d_{j,k}(\xi) &= 2^{km_j}\omegaF_k(\xi),
\end{aligned}
\end{equation}
for $j=3,\ldots,N$.
\\

Using these operators one can show \cite[page 26]{RRS2} that for any $M\geq 0$, the operator $T_{\sigma_{1,2}}^{\Phi}$ can be written as
\begin{multline*}\label{sigma12 original}
T_{\sigma_{1,2}}^{\Phi}(f_1, \ldots, f_N)\\ =\int \sum_{k\geq k_0}^\infty  M_{\mathfrak m}\circ T^{\phase_0}_{d_0} \circ P_{k}^0\left[(Q_{k}^{u_1} \circ T^{\phase_1}_{d_1})(f_1)\, (Q_{k}^{u_2} \circ T^{\phase_2}_{d_2})(f_2)\,\prod_{j=3}^N(P_{k}^{u_j} \circ T^{\phase_j}_{d_{j,k}})(f_j)\right]\\ \times \frac1{(1+|U|^2)^M}  \, \ddd U,
\end{multline*}
for a certain smooth function $\mathfrak{m}$  depending on $\sigma_{1,2},$ with uniformly bounded derivatives of all orders.
Therefore one can reduce the analysis of boundedness of $T_{\sigma_{1,2}}^\Phase$, to the study of the boundedness of the multilinear operator 
\nma{t_piecedecomp for sigma N+1}{
 & D(f_1,\ldots,f_N)(x) \\
 &=\sum_{k\geq k_0}^\infty  M_{\mathfrak m}\circ T^{\phase_0}_{d_0} \circ P_{k}^0\left[(Q_{k}^{u_1} \circ T^{\phase_1}_{d_1})(f_1)\, (Q_{k}^{u_2} \circ T^{\phase_2}_{d_2})(f_2)\,\prod_{j=3}^N(P_{k}^{u_j} \circ T^{\phase_j}_{d_{j,k}})(f_j)\right](x),
}
see \cite{RRS2} for further details.

\section{A catalogue of end-point cases}\label{endpoint cases}

The method by which we prove the boundedness of the components $T_{\sigma_0}^\Phase$, $T_{\sigma_1}^\Phase$ and $T_{\sigma_{1,2}}^\Phase$ splits into four separate cases. For $T_{\sigma_1}^\Phase$ and $T_{\sigma_{1,2}}^\Phase$ we use vector-valued inequality techniques to deal with almost all function spaces $X^p$. However, as mentioned in the introduction, this method fails when $p=2$ or $p=\infty$, so we make use of different techniques when these functions spaces are present. This failure is due in the first case to a lack of usable decay in the amplitude and in the second case due to a lack of a suitable characterisation of $\bmo$, and means we use three different techniques to deal with $T_{\sigma_1}^\Phase$ and $T_{\sigma_{1,2}}^\Phase$. Finally, we make use of a fourth method, which deals with $T_{\sigma_0}^\Phase$ for all values of the function space exponents.

As far as boundedness of $T_{\sigma_1}^\Phase$ is concerned, due to the symmetry of \eqref{the magic operator} in the indices $j=2,\dots,N,$ as was shown in \cite{RRS2}, we only need to consider endpoint cases $(p_0,\dots,p_N)$ which are distinct within the equivalence class of permutations of $(p_2,\dots,p_N)$. Thus, there are three possibilities for the function space with exponent $p_1$: $h^{p_1}$, $L^2$ or $\bmo$. Then for the exponents $p_2,\dots,p_N$ we can have a Cartesian product of the same spaces:
\begin{equation*}
    \left(\prod_{j\in\mathcal{I}_2} L^2\right) \times  \left(\prod_{j\in\mathcal{I}_\infty} \bmo\right) \times
    \left(\prod_{j\in\mathcal{I}_f} h^{p_j}\right),
\end{equation*}
where the index sets $\mathcal{I}_2$, $\mathcal{I}_\infty$ and $\mathcal{I}_f$ are the sets of all $j$ such that $p_j=2$, $p_j=\infty$ and $p_j$ is any other value, respectively.

Similarly, regarding $T_{\sigma_{1,2}}^\Phase$, due to the symmetry of the form of \eqref{t_piecedecomp for sigma N+1} in the indices $j=1,2$ and $j=3,\dots,N$, we only need to consider endpoint cases $(p_0,\dots,p_N)$ which are distinct within the equivalence class of permutations of $(p_1,p_2)$ and $(p_3,\dots,p_N)$. (We have, therefore, $(3+3) \times (3+3+1) - 3 = 39$ cases, since the possibility of three or more copies of $L^2$ appearing is ruled out because $p_0 > 2/3$.)\\

For the Banach target spaces, both for later use in Section \ref{space-time} and for the convenience of the reader, we recall the endpoint-cases that need to be considered and the orders of decay of the amplitude that are involved in each case.  This is of course quite similar to the analysis that was carried out in \cite[Section 5]{RRS2}, with the only difference that here we also consider the cases of various multilinear OIOs. However the interpolation procedure towards the establishment of Banach-target results remain the same. We summarise this in the following lemma:

\begin{lem}\label{cor:endpointcases}

Let $m=\sum_{j=0}^{N} m_j,$ $\frac{1}{p_0}= \sum_{j=1}^{N} \frac{1}{p_j}$, and $\sigma(x,\Xi)\in S^m(n,N)$ and $\varphi_j$ be phase functions of order $s$ with $s>0.$ Let also 
\begin{equation}\label{m}
    m(p)=\begin{cases}
     -(n-1)\abs{\frac{1}{p}-\frac{1}{2}}, \, n>1, \,\varphi_j \mathrm{'s\, \, positively \,\, homogeneous\,\, of\,\, degree \,\,} 1,\,\, \frac{n}{n+1}<p<\infty \\
     -ns\abs{\frac{1}{p}-\frac{1}{2}},\,  \,\varphi_j \mathrm{'s\, \, of \,\, order \,\,} s,\,\, \frac{n}{n+\min(1,s)}<p<\infty.
    \end{cases}
\end{equation}
For Banach-target spaces $($i.e. $X^p$ with $p\in [1,\infty]$$)$, it is enough to prove \emph{Theorem \ref{thm:main}} for the following values of exponents:
\begin{enumerate}[label=\emph{(}\roman*\emph{)}, ref=(\roman*)]
    \item \label{test1} {\bf{Target $\bmo$.}} 
$    \prod_{j=1}^{N} \bmo\to \bmo,$ i.e.
    $(p_j, m_j)=(\infty,m(\infty))$  for all $j=1,\ldots N;$\\
       
    \item {\bf{Target $L^2$.}}  $(p_0,m_0)=(2,0)$, and for each $1\leq j\leq N$, $(p_{j}, m_j)=(2,0)$ and $(p_k, m_k)=(\infty,m_k(\infty))$ for $k\neq j;$\\
        \item {\bf{Target $h^1$.}}  $(p_0, m_0)= (1, m(1))$ and any pair $1\leq j_1< j_2\leq N$,  $(p_{j_1}, m_{j_1})=(p_{j_2}, m_{j_2})=(2,0)$ and $(p_k, m_k)=(\infty, m(\infty))$ for $j_1\neq k\neq j_2$; and
    \item {\bf{Target $h^1$.}} $(p_0, m_0)= (1, m(1))$ and for any $1\leq j\leq N$, $(p_{j}, m_j)=(1, m(1))$, and $(p_k, m_k)=(\infty, m(\infty))$ for $k\neq j$.

\end{enumerate}
\end{lem}

\begin{proof}

This is a standard application of multilinear interpolation, as was also done in \cite{RRS2}. In short, we take two end points, $P_A= (p_{A,1},\ldots,p_{A,N})$ and $P_B$, from the list above, with corresponding amplitude orders $m_A = ∑_{j=1}^N m(p_{A,j})$ and likewise for $m_B$. We then form the amplitude family $σ_z$ given by $σ_z(x,Ξ) = \mathring{σ}(x,Ξ)\jap{Ξ}^{(1-z)m_A + zm_B}$, where $\mathring{σ}\in S^0_{1,0}(n,N)$ is arbitrary, so that $σ_0 \in S^{m_A}_{1,0}$ and $σ_1 \in S^{m_B}_{1,0}$. Notice that for any Schwartz $f_1,\ldots,f_N$, the map $z \mapsto T^Φ_{σ_z}(f_1,\ldots,f_N)$ is analytic, and that the bounds in our proof depend polynomially on $\Im z$. This ensures that we can use the mentioned interpolation result, showing the boundedness of $T^Φ_{σ_z}$ for $z \in [0,1]$. Since $\mathring{σ}$ was arbitrary boundedness holds for any $T^Φ_{σ_z}$ with amplitude $σ_z \in S^{(1-z)m_A + zm_B}_{1,0}(n,N)$, $z \in [0,1]$ and source space $X^{p_1}×\cdots×X^{p_N}$ where $P^{-1} = (p_1^{-1},\ldots,p_N^{-1}) = (1-z)P_A^{-1} + zP_B^{-1}$. One can then do this for any two points $P_A$ and $P_B$ in the convex polygon of studied $P$ to get the full range of exponents, but it suffices to show boundedness at corners and where the function $m$ ceases to be linear. 
\end{proof}

\section{Boundedness of the multilinear operators }\label{sec main multilinear}
In this section we shall very briefly indicate the modifications that are needed in the proofs that were provided in \cite{RRS2}, in order to prove the corresponding results for multilinear OIOs.\\

As far as boundedness results are concerned, due to the symmetry of \eqref{the magic operator} in the indicies $j=2,\dots,N,$ as was shown in \cite{RRS2}, we only need to consider endpoint cases $(p_0,\dots,p_N)$ which are distinct within the equivalence class of permutations of $(p_2,\dots,p_N)$. Similarly, due to the symmetry of the form of \eqref{t_piecedecomp for sigma N+1} in the indicies $j=1,2$ and $j=3,\dots,N$ we only need to consider endpoint cases $(p_0,\dots,p_N)$ which are distinct within the equivalence class of permutations of $(p_1,p_2)$ and $(p_3,\dots,p_N)$.

This reduces the analysis of boundedness of $T^{\Phi}_{\sigma}$ to the investigation of just one of the $T^{\Phi}_{\sigma_j}$)'s say $T^{\Phi}_{\sigma_1},$ one of $T^{\Phi}_{\sigma_{j,k}}$'s say $T^{\Phi}_{\sigma_{1,2}}$ and of course also the boundedness of low-frequency part $T_{\sigma_{0}}^\Phase.$
All the other cases can be studied in essentially identical ways as these.\\
In each case we fix 
\[
    \frac{1}{p_0}=\sum_{j=1}^{N }\frac{1}{p_j},\qquad 1\leq p_j\leq \infty,\quad j=0,\ldots,N,
\]
and $ m_j:=m(p_j)$,  $j=0,\ldots, N,$ with $m(p_j)$ given as in \eqref{m} and consider $f_j \in X^{p_j}$ for $j=1,\dots,N$.
 The rest of the analysis is identical to that of multilinear FIOs as carried out in Section 8 of \cite{RRS}, 
Having this lemma at our disposal, we can run the machinery of the proofs in the case of multilinear FIOs and obtain the desired results.\\

\subsection{Boundedness of $T_{\sigma_{0}}^\Phi$}\label{boundedness of Tsigma0}
Here, due to the localised nature of the amplitude and in contrast to the other parts of the OIO, we can furnish a proof which covers both the quasi-Banach and Banach target spaces cases. In order to control $T_{\sigma_{0}}^\Phi$ defined in \eqref{final form of sigma0}, we observe that since $\theta\in\mathcal{C}_c^{\infty}(\R^n)$, Lemma \ref{low freq lemma llf 1} yields that
\begin{equation*}
    \norm{T_{\theta}^{\phase_j}(f)}_{X^p} \lesssim \norm{f}_{X^p} \quad \mbox{and} \quad \norm{T_{\theta(\cdot/\sqrt{N})}^{\phase_0}(f)}_{X^p} \lesssim \norm{f}_{X^p}
\end{equation*}
for $n/(n+s_c)<p\leq \infty$. Applying these two estimates, the fact that each term is frequency localised, the translation invariance of the norms and H\"older's inequality (using the Littlewood--Paley characterisation of local Hardy spaces) altogether yield
\[
    \norm{\left(\prod_{j=1}^N T_{\theta}^{\phase_j}\circ\tau_{\frac{2\pi k_j}{L}}(f_j)\right)}_{h^r}\lesssim \prod_{j=1}^N \norm{f_j}_{h^{p_j}}.
\]
Combining these estimates one has 
\[
    \norm{T_{\theta(\cdot/\sqrt{N})}^{\phase_0}\left(\prod_{j=1}^N T_{\theta}^{\phase_j}\circ\tau_{\frac{2\pi k_j}{L}}(f_j)\right)}_{X^{p_0}}\lesssim   \prod_{j=1}^N \norm{f_j}_{X^{p_j}}, 
\]
for all the endpoint cases of $p_0,p_1,\dots,p_N$ in Lemma~\ref{cor:endpointcases}. Finally, the boundedness of $T_{\sigma_0}^{\phase_0}$ follows by applying \eqref{derivative of a_K} with the inclusions $\mathcal{C}^1_b\cdot h^p\subseteq h^p$, $L^\infty \cdot L^2 \subseteq L^2$ and $\mathcal{C}^1_b\cdot \bmo \subseteq \bmo$ (see \cite{Gold}).

Therefore, for the purely low-frequency portion of the operator, we have now established the boundedness with both Banach and quasi-Banach target spaces.

\subsection{Boundedness of $T_{\sigma_{1}}^\Phase$}
Due to the symmetry of the representation \eqref{the magic operator} of $T_{\sigma_{1}}^\Phase$ (in the indicies $j=2,\dots,N$) we only need to consider endpoint cases $(p_0,\dots,p_N)$ which are distinct within the equivalence class of permutations of $(p_2,\dots,p_N)$.

\subsubsection{Boundedness with Banach targets.}
Having this, then all the boundedness results with target spaces $L^2$ and $\bmo$ (in accordance to Theorem \ref{cor:endpointcases}) are proven in exactly the same way as in the case of multilinear FIOs in \cite{RRS2}, where one replaces $-(n-1)/2$ of multilinear FIOs by $-ns/2$ of multilinear OIOs and noting that no restriction on the dimension (as in the FIO case) is necessary, since $-ns/2<0$. 

\subsubsection{Boundedness with quasi-Banach targets.}
As discussed earlier in connection to representation \eqref{the magic operator}, matters can be reduced to the study of the regularity of the multilinear operator
\begin{equation}\label{defn:thick I}
\begin{split}
\textbf{I}:=  \sum_{k=k_0}^{\infty} \,Q_{k}^0 \,\left[ (Q_{k}^{u_1} \circ T^{\phase_1}_{b_1})(f_1)\,\prod_{j=2}^{N} (P_{k}^{u_j} \circ T^{\phase_j}_{b_{j,k}})(f_j)\right].
\end{split}
\end{equation}

Our goal is to prove the boundedness of $T^\Phi_{\sigma_1}$ with target in $h^{p_0}$ with $n/(n+ s_c)<p_0<\infty$ and $p_0\neq 2$. We also note that the cases $p_0\geq 1$ are all Banach, but our method of proof will cover these cases as well.  Using \eqref{defn:thick I}, we infer that the boundedness of $T^\Phi_{\sigma_1}$, could via \eqref{norm in loc hardy}, be investigated by considering
\begin{equation}\label{eq:low_freq_theta_2}
\begin{split}
    \vartheta_0(D)(\textbf{I})=\sum_{k=k_0}^{3}   \vartheta_0(D)\, Q_k^0(D)\left[ (Q_{k}^{u_1} \circ T^{\phase_1}_{b_1})(f_1)\, \prod_{\ell=2}^{N} (P_{k}^{u_\ell} \circ T^{\phase_\ell}_{b_{\ell,k}})(f_\ell)\right],
    \end{split}
 \end{equation}
and for $j\geq 1$
\begin{equation*}\label{eq:high_freq_theta_2}
    \begin{split}
    \vartheta_j(D)(
    \textbf{I})&= \sum_{k=k_0,|k-j|\leq 4}^{\infty}  \vartheta_j(D)\, Q_k^0(D) \left[ (Q_{k}^{u_1} \circ T^{\phase_1}_{b_1})(f_1)\, \prod_{\ell=2}^{N} (P_{k}^{u_\ell} \circ T^{\phase_\ell}_{b_{\ell,k}})(f_\ell)\right]\\
    &=\sum_{\ell=-4}^4
     \mathbb{1}_{[k_0,\infty)}(\ell+j) \vartheta(2^{-j}D)\, {\phi}(2^{-(j+\ell)}D) [F_{j+\ell}^U],
    \end{split}
\end{equation*}
where for all $k\in \Z$ 
\[
    {F_k^U}=(Q_{k}^{u_1} \circ T^{\phase_1}_{b_1})(f_1)\, \prod_{\ell=2}^{N} (P_{k}^{u_\ell} \circ T^{\phase_\ell}_{b_{\ell,k}})(f_\ell).
\]
Now given an $(N-1)$-tuple $(p_2,\ldots, p_N),$ we define 
\[
    \mathfrak{I}_2=\{j\in \{2,\ldots,N\}:\, p_j=2\},\, \mathfrak{I}_f=\{j\in \{2,\ldots,N\}:\, 2\neq p_j<\infty\}
\]
and 
\[
    \mathfrak{I}_\infty=\{j\in \{2,\ldots,N\}:\, p_j=\infty\}.
\]
Using this notation we can write
\begin{equation}\label{defn:fuk}
    {F_k^U}:=(Q_{k}^{u_1} \circ T^{\phase_1}_{b_1})(f_1)\,\prod_{j\in \mathfrak{I}_2} (P_{k}^{u_j} \circ T^{\phase_j}_{b_{j,k}})(f_j) 
    \prod_{j\in \mathfrak{I}_f} (P_{k}^{u_j} \circ T^{\phase_j}_{b_{j,k}})(f_j)\prod_{j\in \mathfrak{I}_\infty} (P_{k}^{u_j} \circ T^{\phase_j}_{b_{j,k}})(f_j).
\end{equation}

Taking \eqref{norm in loc hardy} into account for a generic piece of $F^{U}_k$, we shall see that the following proposition will be useful in dealing with various cases that arise in connection to the proof of $h^{p_0}$-regularity of $\bf{I}$ (given by \eqref{defn:thick I}).
\begin{prop}\label{prop:Tedesco}
 Given $p_1>0$ and $N_1\geq 2$, assume that $\frac{1}{r_0}= \frac{1}{p_1}+ \frac{N_1-1}{2}$ and $X^{p_1}$ be defined as in \eqref{defn:xp}. Then one has\[
       \Vert\Big(\sum_{k=k_0}^{\infty}\,  |(Q_k^{u_1}f_1) \prod_{j=2}^{N_1} (P_{k}^{u_j} (f_j)|^2\Big)^{1/2}\Big\Vert_{L^{r_0}}\lesssim \norm{f_1}_{X^{p_1}}\prod_{j=2}^{N_1} \norm{f_j}_{2}.
\]
\end{prop}
\begin{proof} 
By the translation invariance of the norm of the spaces $X^p$, we can reduce the study to the case where $u_1=\ldots=u_{N_1}=0$.

Consider the multilinear pseudodifferential operator
 \[
  T_k\left(f_{1},\ldots, f_{N_1}\right)(x):= Q_{k} f_1\,\prod_{j=2}^{N_1} P_{k} f_j
 \]
 with the symbol 
 \[
    \rho^{k}(\Xi)=\psi(2^{-k}\xi_1)\,\,\prod_{j=2}^{N_1}\theta(2^{-k}\xi_j).
 \]
 
 Now since, in addition to the frequency localisations $|2^{-k} \xi_1|\sim 1$, $|2^{-k} \xi_j|\lesssim 1$ for $2\leq j\leq N_1$, one also has that $|\Xi|\leq c|\xi_1|$ on the support of $\sigma_1$ (note that the later follows from \eqref{biglambda} and \eqref{defn:sigmaj}), then the Leibniz rule, the aforementioned support properties, and finally \eqref{the key on the support of sigma1} yield
\begin{equation*}\label{andalucia2}
\begin{split}
    \vert &\partial^{\alpha}_{\Xi}(\rho^k(\Xi)) \vert\lesssim \langle \Xi\rangle^{-|\alpha|},
\end{split}
\end{equation*}
 which yields that $\rho^k\in S^0_{1,0}(n,N),$ uniformly in $k$.
 
Let assume first that $p_1<\infty$.  Khinchin's inequality yields that
\begin{equation*}\label{rascal khinchin}
\begin{split}
    \Big\Vert\Big(
        \sum_{k=-5}^{\infty}\,  |(Q_k f_1) \prod_{j=2}^{N_1} (P_{k} f_j)|^2\Big)^{1/2}\Big\Vert_{L^{r_0}}&\approx
\norm{\sum_{k\geq -5}
	  \varepsilon_k (t)
	  Q_k f_1\prod_{j=2}^{N_1} P_k f_j}_{L_{x,t}^{r_0}(\R^n\times [0,1])} \\&= \norm{\sum_{k\geq -5}
	  \varepsilon_k (t)
	 T_{k}(f_1, \dots , f_{N_1})}_{L_{x,t}^{r_0}(\R^n\times [0,1])},
	 \end{split}
\end{equation*}
where $\{\varepsilon_j (t)\}_j$ are the Rademacher functions.
Now the family of multilinear pseudodifferential operators 
$\sum_{k=-5}^{\infty}  \varepsilon_k(t) T_{k}(f_1, \dots , f_{N_1}),$ has the symbol
\begin{equation*}\label{jaevla rho}
    \rho_{t}         (\xi_1,\ldots,\xi_{N_1}):=\sum_{k=-5}^\infty \varepsilon_k(t) \rho^k(\xi_1,\dots,\xi_{N_1}) \in S^0_{1,0}(n,N),
\end{equation*}
uniformly in $t$. Therefore, the boundedness of multilinear pseudodifferential operators of order zero from $\prod_{j=1}^{N} h^{l_j}\to L^r$  \cite{TanZhao}*{Theorem 1.1} yields
\[
    \norm{\sum_{k\geq -5}
	\varepsilon_k (t) T_{k}(f_1, \dots , f_{N_1})}_{L_{x,t}^{r_0}(\R^n\times [0,1])}\lesssim \norm{f_1}_{h^{p_1}} \prod_{j=2}^{N_1} \norm{f_j}_{L^2}.
\]

Now let us assume now that $p_1=\infty$. Note that we are also allowed to assume \eqref{ineqforcut-off} on the support of $\rho^{k}$, which yields that 
\begin{equation*}
\begin{split} 
&T_k(f_1,\ldots,f_{N_1})(x)\\
&=
\int_{\left(\mathbb{R}^{nN_1}\right)}  2^{N_1 n k}\phi^{\vee}\left(2^{k}\left(x-y_{1}\right), \ldots, 2^{k}\left(x-y_{N_1}\right)\right) \prod_{j=2}^{N_1} P_{k}^{}\left(f_{j}\right)\left(y_{j}\right) Q_{k}^{}\left(f_{1}\right)\left(y_{1}\right) d Y.
\end{split}
\end{equation*}
H\"older's inequality, the translation invariance of the Lebesgue measure and the definition of the maximal operator $\mathfrak{M}^{p}_{a,b}$ in \eqref{Parkmax} yield
\begin{equation*}\label{pointwisebloodymess}
    \begin{split}
    &\abs{T_k(f_1,\ldots, f_{N_1})(x)}\leq\\\nonumber
    &2^{N_1nk}\left\{\int_{\R^{nN_1}}\jap{2^k Y}^{sq'} \abs{\phi^{\vee}(2^{k}Y)}^q\dd Y \right\}^{1/q'}
    \left\{\int_{\R^{nN_1}} \prod_{j=2}^{N_1} \frac{\abs{P_{k}^{}\left(f_{j}\right)\left(y_{j}\right)}^q}{\jap{2^k(x-y_j)}^{sq/N_1}}
    \frac{\abs{Q_{k}^{}\left(f_{1}\right)\left(y_{1}\right)}^q}{\jap{2^k(x-y_1)}^{sq/N_1}}\dd Y
    \right\}^{1/q}\\\nonumber
    &\lesssim 2^{-N_1nk/q'} 2^{N_1nk}\mathfrak{M}_{s/N_1,2^k}^q(Q_k f_1)(x) 2^{-k n/q} \prod_{j=2}^{N_1} \mathfrak{M}_{s/N_1,2^k}^q (P_k f_j)(x) 2^{-k n/q} \\
    \nonumber&\lesssim 
    \mathfrak{M}_{s/N_1,2^k}^q(Q_k f_1)(x) \prod_{j=2}^{N_1} \mathfrak{M}_{s/N_1,2^k}^q (P_k f_j)(x)
\end{split}
\end{equation*}
where we have also used that for all $z\in \R^{nN_1}$
\[
    (1+2^{2\ell}\abs{z_1}^2+\ldots+2^{2\ell}\abs{z_{N_1}}^2)^{N_1}\geq  \prod_{k=1}^{N_1}(1+2^{2\ell}\abs{z_k}^2).
\]
Now denoting the set of all dyadic cubes in $\R^n$ by $\mathcal{D}$, and denoting for each $k\in \Z$ the elements of $\mathcal{D}$ with side length $2^{-k}$ by $\mathcal{D}_k,$ we have by inequality \eqref{Parks inequality} that
for every dyadic cube $J \in \mathcal{D}_k$ and every $f$ 
\[
    \sup_{y\in J} \mathfrak{M}_{s/N,2^k}^q(f)(y)\lesssim \inf_{y\in J} \mathfrak{M}_{s/N,2^k}^q(f)(y),
\]
with constants independent of $f$ and $k$. 

Therefore, since there is no overlap between $\mathcal{D}_k$'s, we have
\begin{equation}\label{eq:Park}
    \begin{split}
     \Big\Vert\Big(\sum_{k=-5}^{\infty}\,  &|T_k(f_1,\ldots f_{N_1})|^2\Big)^{1/2}\Big\Vert_{L^{r_0}}=\Big\Vert\Big(
\sum_{k=-5}^{\infty}\,  \sum_{J\in \mathcal{D}_k}|T_k(f_1,\ldots f_{N_1})|^2\chi_J \Big)^{1/2}\Big\Vert_{L^{r_0}} \\
&\leq \norm{\left(\sum_{k\geq -5}\sum_{J\in \mathcal{D}_k} \prod_{j=2}^{N_1}  \abs{\mathfrak{M}_{s/N_1,2^k}^q (P_k f_j)(x)}^2\abs{\mathfrak{M}_{s/N_1,2^k}^q(Q_k f_1)}^2\chi_J\right)^{1/2}}_{L^{r_0}}\\
&\leq \norm{\left(\sum_{k\geq -5}\sum_{J\in \mathcal{D}_k} \prod_{j=2}^{N_1}  \inf_{y\in J}\abs{\mathfrak{M}_{s/N_1,2^k}^q (P_k f_j)(x)}^2\abs{\inf_{y\in J}\mathfrak{M}_{s/N_1,2^k}^q(Q_k f_1)}^2\chi_J\right)^{1/2}}_{L^{r_0}}\\
&\leq \norm{\left(\sum_{k\geq -5} \brkt{\sum_{J\in \brkt{\mathcal{D}_k}} \prod_{j=2}^{N_1}  \inf_{y\in J}\mathfrak{M}_{s/N_1,2^k}^q (P_k f_j)(x)\inf_{y\in J}\mathfrak{M}_{s/N_1,2^k}^q(Q_k f_1)\chi_J}^2\right)^{1/2}}_{L^{r_0}}.
    \end{split}
\end{equation}

Now by Theorem \ref{Parks thm},  given $r_0\in (0,\infty]$, $0<q\leq \infty$, $\gamma\in (0,1)$ and $s/N_1 >n/(\min\brkt{2,q, r_0}),$ for any dyadic cube $J\in \mathcal{D}$ there exists a measurable subset $S_J\subset Q$, depending on $\gamma, f_k ,q,s,N_1$ such that $\abs{S_J}>\gamma \abs{J}$. For this $S_J$ and any $0<\rho<\infty$ one has for \(x \in J\) that
\[
    \chi_J(x)=1<\frac{1}{\gamma^{1 / \rho}} \frac{\left|S_{J}\right|^{1 / \rho}}{|J|^{1 / \rho}}=\frac{1}{\gamma^{1 / \rho}}\left(\frac{1}{|J|} \int_{J} \chi_{S_{Q}}^{\rho}(y) d y\right)^{1 / \rho} \leq \gamma^{-1 / \rho} \mathcal{M}_{\rho}\left(\chi_{S_{J}}\right)(x).
\]

Hence, using this and the vector-valued maximal inequality \eqref{FSvector}, one can bound the last term in \eqref{eq:Park} by 
\begin{equation}\label{fanskapet}
 \norm{\left(\sum_{k\geq -5} \brkt{\sum_{J\in {\mathcal{D}_k}} \prod_{j=2}^{N_1}  \inf_{y\in J}\mathfrak{M}_{s/N_1,2^k}^q (P_k f_k)(x)\inf_{y\in J}\mathfrak{M}_{s/N_1,2^k}^q(Q_k f_1)\chi_{S_J}}^2\right)^{1/2}}_{L^{r_0}}.   
\end{equation}

We also note that the characterisation of BMO given in Theorem \ref{Parks thm} implies that, given $0<q\leq \infty$, $\gamma\in (0,1)$ and $s/N_1 >n/(\min\brkt{2,q}),$ one has 
\[
    \norm{\left(\sum_{k\geq -5} \brkt{\sum_{J\in {\mathcal{D}_k}} \inf_{y\in J}\mathfrak{M}_{s/N_1,2k}^q(Q_k f_1)\chi_{S_J}}^2\right)^{1/2}}_{L^{\infty}}\approx \norm{\Gamma(D)f_1}_{\mathrm{BMO}}\leq \norm{f_1}_{\bmo}.
\]
where $\Gamma(D)$ is a high-frequency cut-off.

Therefore, H\"older's inequality, Theorem \ref{Park2} and the $L^2$-boundedness of Hardy-Littlewood's maximal functions yield that the expression in \eqref{fanskapet} is bounded by
\[
    \begin{split}
    &\norm{\prod_{j=2}^{N_1}\sup_k \mathfrak{M}_{s/N_1,2^k}^q (P_k f_j)}_{L^{r_0}} \norm{\left(\sum_{k\geq -5} \brkt{\sum_{J\in \brkt{\mathcal{D}_k}} \inf_{y\in J}\mathfrak{M}_{s/N_1,2^k}^q(Q_k f_1)\chi_{S_J}}^2\right)^{1/2}}_{L^{\infty}}\\
    &\lesssim\prod_{j=2}^{N_1}\norm{\sup_k \mathfrak{M}_{s/N_1,2^k}^q (P_k f_j)}_{L^{2}} \norm{f_1}_{\mathrm{bmo}}\lesssim \prod_{j=2}^{N_1}\norm{\sup_k |P_k f_j |}_{L^{2}} \norm{f_1}_{\mathrm{bmo}}\\&\lesssim \prod_{j=2}^{N_1}\norm{\mathcal{M}{f_j}}_{L^{2}}\norm{f_1}_{\bmo}\lesssim \prod_{j=2}^{N_1}\norm{f_j}_{L^{2}}\norm{f_1}_{\bmo}.
    \end{split}
\]

\end{proof}
Now we turn to the study of the regularity of the multilinear operators associated to $T_{\sigma_1}^{\Phi}$. This will be divided in the following cases:\\

{\bf{Case I. $\mathfrak{I}_2\neq \emptyset $.}} Observe that by our previous considerations the frequency support of  $F_k^U$ (given by \eqref{defn:fuk}) is contained in  $B(0, 2^k R)$, for some $R\geq 1.$ Therefore, for $\ell\in[-4,4]$, the frequency support of $F_{j+\ell}^U$ is contained in $B(0,2^j (2^{\ell}R))$. 
Hence using \eqref{t_Triebels vector valued} we have
\begin{equation*}\label{t_?}
    \Big\Vert \Big(\sum_{j=1}^{\infty}|\vartheta_j(D)(\textbf{I})|^2\Big)^{1/2}\big\Vert_{L^{p_0}} \lesssim
     \Big\Vert\Big(
    \sum_{k=-5}^{\infty}\,  |F^{U}_{k}|^2\Big)^{1/2}\Big\Vert_{L^{p_0}}.
\end{equation*}
 Note that for $j\in \mathfrak{I}_2$, the $b_{j,k}$'s dependence on $k$ could be suppressed due to the fact that for these terms the corresponding $m_j$'s are equal to zero and one can replace the amplitudes by the constant function one.
 Hence using the uniform boundedness given in  \eqref{eq:uniform_bounddness}, the embedding $\ell^1(\N)\subset \ell^2(\N)$ jointly with the Cauchy--Schwarz inequality, H\"older's inequality, Lemma \ref{Lu-Be}, Proposition \ref{prop:Tedesco} and the boundedness of linear oscillatory integrals given in Theorem \ref{linearhpthmoio},  yield
 \begin{equation}\label{eq:FingMess}
    \begin{split}
     &\Big\Vert\Big(\sum_{k=-5}^{\infty}\,  |F_k^U|^2\Big)^{1/2}\Big\Vert_{L^{p_0}}\\
    &\lesssim
    \Big\Vert\Big(
    \sum_{k=-5}^{\infty}\,  |(Q_k^{u_1}T_{b_1}^{\varphi_1})f_1) \prod_{\mathfrak{I}_2} (P_{k}^{u_j} T^{\phase_j}_{1})(f_j)|^2\Big)^{1/2}\Big\Vert_{L^{r_0}}\\
        &\qquad\times\prod_{\mathfrak{I}_f}\Big\Vert\Big(\sum_{k=-5}^{\infty}\,  |(P_k^{u_j} T_{b_{j,k}}^{\varphi_j})(f_j)|^2\Big)^{1/2}\Big\Vert_{L^{p_l}}  \prod_{ \mathfrak{I}_\infty} \norm{f_j}_{\bmo}\\
    &\lesssim \norm{T_{b_1}^{\varphi_1} f_1}_{X^{p_1}} \prod_{\mathfrak{I}_2} \norm{T^{\phase_j}_{1}(f_j)}_{L^{2}} \prod_{\mathfrak{I}_f} \norm{f_j}_{h^{p_j}} \prod_{\mathfrak{I}_\infty} \norm{f_j}_{\bmo}\lesssim \prod_{j=1}^N \norm{f_j}_{X^{p_j}},
\end{split}
\end{equation}
where
\[
    \frac{1}{p_0}=\frac{1}{r_0}+\sum_{\ell\in \mathfrak{I}_f} \frac{1}{p_\ell}, \qquad \frac{1}{r_0}=\frac{1}{p_1}+\frac{\abs{\mathfrak{I}_2}}{2}.
\]

{\bf{Case II. $\mathfrak{I}_2=\emptyset$.}}
In this case $r_0=p_1$ and if moreover $p_1<\infty$ the we have
\begin{equation}\label{eq:Finite_p}
    \Big\Vert\Big(
    \sum_{k=-1}^{\infty}\,  |(Q_k^{u_1}T_{b_1}^{\varphi_1}) f_1) |^2\Big)^{1/2}\Big\Vert_{L^{r_0}}\lesssim \norm{f_1}_{h^{p_1}},
\end{equation}

and we can proceed as  in \eqref{eq:FingMess} to reach the desired estimate. However, if $p_1=\infty$ then  the classical Fefferman--Stein estimate yields that
\begin{equation}\label{basic feffstein}
  \sup_{k} \norm{Q_k^{u_1}  (f_1)}_{\infty}\lesssim \norm{f_1}_{\bmo}.  
\end{equation}
Hence Lemma \ref{Lu-Be} yields
\begin{equation*}
    \begin{split}
     &\Big\Vert\Big(
\sum_{k=-5}^{\infty}\,  |F_k^U|^2\Big)^{1/2}\Big\Vert_{L^{p_0}}\\
&\lesssim \norm{f_1}_{\bmo}
\prod_{\mathfrak{I}_f}\Big\Vert\Big(
\sum_{k=-5}^{\infty}\,  |(P_k^{u_j} T_{b_{j,k}}^{\varphi_j})(f_j)|^2\Big)^{1/2}\Big\Vert_{L^{p_j}}  \prod_{\mathfrak{I}_\infty} \norm{f_j}_{\bmo} \lesssim  \prod_{j=1}^N \norm{f_j}_{X^{p_j}}.
    \end{split}
\end{equation*}

Finally, for the low frequency part \eqref{eq:low_freq_theta_2}, we only need to estimate the $L^{p_0}$ norm of $F^{U}_k$. To that end, we use the following generalised Hölder's inequality
\begin{equation}\label{eq:General_Holder}
    \Vert \prod_{j=1}^{N} f_j\Vert_{L^{p_0}}\lesssim \prod_{\mathfrak{I}_2 \cup \mathfrak{I}_f} \Vert  f_j\Vert_{h^{p_j}} \prod_{\mathfrak{I}_\infty} \Vert f_j\Vert_{L^{\infty}} ,
\end{equation}
where $\frac{1}{p_0}= \sum_{j=1}^{N}\frac{1}{p_j} ,$
which is is a consequence of \cite{TanZhao}*{Theorem 1.1}, together with Lemma \ref{lem:Wase-Lin}, which concludes the proof.

\subsection{Boundedness of $T_{\sigma_{1,2}}^\Phase$}
In the analysis of the boundedness of $T^\Phase_{\sigma_{j,k}}$, the symmetry of the operators form under permutations of the frequency variables allows us to restrict our attention to just one of the $\sigma_{j,k}$, the argument for all the others being identical. For definiteness, we choose to study $\sigma_{1,2}$, so we have that $|\xi_1|$ and $|\xi_2|$ are comparable to each other.

\subsubsection{Boundedness with Banach targets}
The demonstrations of the boundedness of $T_{\sigma_{1,2}}^\Phase$  with target spaces $\bmo$ and $L^2$ are idential to that of multilinear FIOs as carried out in Section 8 of \cite{RRS}. However the analysis in \cite{RRS} required a result about Carleson measures associated to linear FIOs. The analogue of that result was provided in Proposition \ref{lem:smallcarlnorm} above, and with that proposition at hand, we can run the machinery of the proofs in the case of multilinear FIOs in \cite{RRS2} and obtain the boundedness of $T_{\sigma_{1,2}}^\Phase$ with target spaces $\bmo$ and $L^2$.

\subsubsection{Boundedness with quasi-Banach targets.}
Using the representation \eqref{t_piecedecomp for sigma N+1},  we are dealing with the $h^p_0$-regularity of the multilinear operator
\begin{equation}\label{t_piecedecomp for sigma N+1 bilin}
   D(f_1, \ldots,f_N)(x)=\sum_{k=k_0}^\infty  M_{\mathfrak m}\circ T^{\phase_0}_{d_0} \circ P_{k}^0\left[G^{U}_{k}\right](x), 
\end{equation}
 
where 
\[
\begin{split}
    G^{U}_{j}:=(Q_{j}^{u_1} \circ T^{\phase_1}_{d_1})(f_1)(Q_{j}^{u_2} \circ T^{\phase_2}_{d_2})(f_2)\prod_{\iota\in \mathfrak{I}_2}(P_{j}^{u_\iota} \circ T^{\phase_\iota}_{d_{\iota,j}})(f_\iota)\\\times \prod_{\iota\in \mathfrak{I}_f}(P_{j}^{u_\iota} \circ T^{\phase_\iota}_{d_{\iota,j}})(f_\iota) \prod_{\iota\in \mathfrak{I}_\infty}(P_{j}^{u_\iota} \circ T^{\phase_\iota}_{d_{\iota,j}})(f_\iota)
    \end{split}
\]
Now given an $(N-2)$-tuple $(p_3,\ldots, p_N),$ we define 
\[
    \mathfrak{J}_2=\{j\in \{3,\ldots,N\}:\, p_j=2\},\, \mathfrak{J}_f=\{j\in \{3,\ldots,N\}:\, 3\neq p_j<\infty\}
\]
and 
\[
    \mathfrak{J}_\infty=\{j\in \{3,\ldots,N\}:\, p_j=\infty\}.
\]
Now, by using \eqref{decomposition of bj} we can rewrite
\[
    d_0(k,\xi)=d_0^\flat(k,\xi)+d_0^\sharp(k,\xi),
\]
and observe that the supports of $\omega_k$ and $\chi_0$ allow us to write 
\[
    T_{d_0^\flat}^{\phase_0}\circ P_k^0=2^{km_0}T_{\Omega_0}^{\phase_0}\circ P_{k_0}^0,
\]
where $\Omega_0:=1-\chi_0$.

\noindent{\bf The analysis concerning $d_0^\sharp$.}
For $d_0^\sharp(k,\xi)$, we replace $(T^{\phase_0}_{d_0^\sharp}\circ P_{k}^{0})(f)$ by $ T^{\phase_0}_{\gamma} \circ R_k \circ P_k^{0}(f)$, where $\gamma(\xi):= \chi_0(\xi)|\xi|^{m_0} \in S^{m_0}$ with $m_0<0$, and $R_k$ is as in \eqref{defn: Rk}.
This yields that 
\begin{equation}\label{eq:Decomposition_crap}
    \begin{split}
        &\sum_{k=k_0}^\infty  M_{\mathfrak m}\circ T^{\phase_0}_{d_0^\sharp} \circ P_{k}^0\left[G^U_{k}\right](x)\\
        &=\sum_{j\geq k_0} \sum_{k\geq j} 2^{(k-j)m_0} M_{\mathfrak m} T_\gamma^{\phase_0} Q_jP_k^0\left[G^U_{k}\right](x)\\
        &=\sum_{k\geq 0} 2^{km_0}\sum_{j\geq k_0}  M_{\mathfrak m_{j+k}} T_\gamma^{\phase_0}Q_jP_{k+j}^0\left[G^U_{k+j}\right](x).\\
    \end{split}    
\end{equation}
{\begin{rem}
Note that here the fact that $m_0<0$ $($which excludes target-space $L^2$$)$ is crucial in the analysis that follows below. 
 \end{rem}}

Now we observe that one can write 
\m{
M_{{\mathfrak m_{j+k}}} \circ T^{\phase_0}_{\gamma} \circ Q_{j}
= \left(\sum_{j'=j-1}^{j+1} T_{j,j',k}^{U}\right) \circ Q_{j} \\
= \left(\sum_{\ell=-1}^1\sum_{j'-j\equiv\ell\Mod{3}} T_{j'+\ell,j',k}^{U}\right) \circ Q_{j}
}
where $T_{j,j',k}^{U}$ is the oscillatory integral with amplitude ${\mathfrak m}(j+k,x,U)\, \gamma (\xi)\,\phiF_{j'}(\xi)$ and phase $\phase_0$. Observe that
\[
\mathcal{T}_{j,k}^U:=\sum_{\ell=-1}^{1}\sum_{j'-j\equiv\ell\Mod{3}} T_{j'+\ell,j',k}^{U}
\]
is periodic in $j$ with period $3$, and is an oscillatory integral with amplitude in $S^{m_0}$ uniformly in $k$.

Thus \eqref{eq:Decomposition_crap} can be rewritten as
\[\label{periodicity again}
    \sum_k 2^{km_0}\sum_{\ell=0}^2 \mathcal{T}_{\ell,k}^{U} \left({D}_{\ell,k}(f_1,\ldots,f_N)\right)(x),
\]
where
\[
{D}_{\ell,k}(f_1,\ldots,f_N)(x) 
    := \chi_0(2D){ \sum_{j\equiv \ell\Mod{3},\,j\geq k_0 } \,Q_{j}^0 P_{k+j}^0 \left[G_{j+k}^U\right](x),}
\]
and $\chi_0$ is the same high-frequency cut-off introduced previously (with a symbol in $S^0$).

For the high-frequency part of the multilinear operator we observe that
\[
  \Vert{D_{\ell,k}}(f_1,\ldots, f_N) \Vert_{h^{p_0}}
  \lesssim \Big\Vert \sum_{j\equiv \ell\Mod{3},\,j\geq k_0 } \,Q_{j}^0 P_{k+j}^0 \left[G^U_{j+k}\right] \Big\Vert_{h^{p_0}}.
\]
Now since the spectrum of $Q_{j}^0 P_{k+j}^0 \left[G^U_{j+k}\right]$ is inside an annulus of size $2^j$, a theorem in Section 2.5.2 on page 79 of \cite{Triebel}, together with  estimate \eqref{t_Triebels vector valued} and finally the Cauchy-Schwarz inequality (using the boundedness of the operators $T^{\phase_1}_{d_1}$ and $T^{\phase_2}_{d_2}$), yield that
\begin{equation}\label{eq:peperoni}
\begin{split}
    &\Big\Vert \sum_{j\equiv \ell\Mod{3},\,j\geq k_0 } \,Q_{j}^0 P_{k+j}^0 \left[G^U_{j+k}\right]\Big\Vert_{h^{p_0}}\\
    &\lesssim \Big\Vert \Big(\sum_{j\equiv \ell\Mod{3},\,j\geq k_0 } \, \Big|Q_{j}^0 P_{k+j}^0 \left[G^U_{j+k}\right]\Big|^2\Big)^{\frac{1}{2}}\Big\Vert_{L^{p_0}}\\
    &\lesssim
    \Big\Vert \Big(\sum_{j\geq k+k_0 } \,  \Big|G^U_{j}\Big|^2\Big)^{\frac{1}{2}}\Big\Vert_{L^{p_0}}.
\end{split}
\end{equation}

Now we proceed by dividing the regularity results into cases which we shall deal with accordingly.\\

\noindent {\bf{Case I.  $p_1<\infty$ }}\\
Let us first assume that $p_2 <\infty$. Here we use the same reasoning as in the paragraph preceding the displayed equation \eqref{eq:FingMess} and note that the left-hand side term of \eqref{eq:peperoni} is bounded by
\begin{equation*}
\begin{split}
    \norm{\brkt{\sum_{j\geq 2k_0 } \,  \Big|(Q_{j}^{u_1} \circ T^{\phase_1}_{d_1})(f_1)\Big|^2}^{1/2}}_{L^{p_1}}\\\times
    \Big\Vert \Big(\sum_{j\geq 2k_0 } \,  \Big|(Q_{j}^{u_2} \circ T^{\phase_2}_{d_2})(f_2)\prod_{\iota \in \mathfrak{J}_2 \cup \mathfrak{J}_f \cup \mathfrak{J}_\infty }(P_{j}^{u_\iota} \circ T^{\phase_\iota}_{d_{\iota,j}})(f_\iota)\Big|^2\Big)^{\frac{1}{2}}\Big\Vert_{L^{r_1}},
    \end{split}
\end{equation*}
where 
\[
    \frac{1}{{p_0}}= \frac{1}{p_1}+\frac{1}{r_1},
\]
and the term $(Q_{j}^{u_2} \circ T^{\phase_2}_{d_2})(f_2)\prod_{\iota \in \mathfrak{J}_2 \cup \mathfrak{J}_f \cup \mathfrak{J}_\infty }(P_{j}^{u_\iota} \circ T^{\phase_\iota}_{d_{\iota,j}})(f_\iota)$ is essentially $F^{U}_k$ given in \eqref{defn:fuk}.
Now since we have that  
\[
     \norm{\brkt{\sum_{j\geq 2k_0 } \,  \Big|(Q_{j}^{u_1} \circ T^{\phase_1}_{d_1}(f_1)\Big|^2}^{1/2}}_{p_1} \lesssim \norm{T^{\phase_1}_{d_1}(f_1)}_{h^{p_1}}\lesssim \norm{f_1}_{h^{p_1}},
\]
the same argument as the one involved in deriving estimate \eqref{eq:Finite_p} for the case $N=2$ and \eqref{eq:FingMess} for $N\geq 3$, yield the desired bound.\\

Now if $p_2= \infty$, then by using Fefferman--Stein's estimate \eqref{basic feffstein}, one extracts the $\Vert f_2\Vert_{\bmo}$ from the left-hand side of \eqref{eq:peperoni} and the remaining term will be $$\Big\Vert \Big(\sum_{j\geq 2k_0 } \,  \Big|(Q_{j}^{u_2} \circ T^{\phase_2}_{d_2})(f_1)\prod_{\iota \in \mathfrak{J}_2 \cup \mathfrak{J}_f \cup \mathfrak{J}_\infty }(P_{j}^{u_\iota} \circ T^{\phase_\iota}_{d_{\iota,j}})(f_\iota)\Big|^2\Big)^{\frac{1}{2}}\Big\Vert_{L^{r_1}},$$
for which the boundedness can be established as was done previously.

\noindent {\bf{Case II. $p_1= \infty$.}}\\
 In this case, applying Fefferman--Stein's estimate \eqref{basic feffstein}, one has that \eqref{eq:peperoni} is bounded by 
\[
    \norm{f_1}_{\bmo}\Big\Vert \Big(\sum_{j\geq 2k_0 } \,  \Big|(Q_{j}^{u_2} \circ T^{\phase_2}_{d_2})(f_2)\prod_{\iota \in \mathfrak{J}_2 \cup \mathfrak{J}_f \cup \mathfrak{J}_\infty }(P_{j}^{u_\iota} \circ T^{\phase_\iota}_{d_{\iota,j}})(f_\iota)\Big|^2\Big)^{\frac{1}{2}}\Big\Vert_{L^{p_0}}.
\]
If we assume that $p_2<\infty$, we just proceed as in the analysis of \eqref{eq:FingMess} (or {\bf{Case I}} above).\\

Now if $p_2 = \infty$, since $p_0<\infty$, then $\mathfrak{J}_2\cup \mathfrak{J}_f\neq \emptyset$. If $\mathfrak{J}_2\neq \emptyset$, \eqref{eq:peperoni} is bounded by
\[
\begin{split}
    \norm{f_1}_\bmo \Big\Vert \Big(\sum_{j\geq 2k_0 } \,  \Big|(Q_{j}^{u_2} \circ T^{\phase_2}_{d_2})(f_2)\prod_{\iota \in \mathfrak{J}_2  }(P_{j}^{u_\iota} \circ T^{\phase_\iota}_{1})(f_\iota)\Big|^2\Big)^{\frac{1}{2}}\Big\Vert_{L^{r_2}} \\\times\prod_{\mathfrak{J}_f}\norm{\Big(\sum_{j\geq 2k_0 } \,  \Big|( P_{j}^{u_\iota} \circ T^{\phase_\iota}_{d_\iota,j})(f_\iota)\Big|^2\Big)^{\frac{1}{2}}}_{L^{p_\iota}} \prod_{\mathfrak{J}_\infty} \norm{f_\iota}_\bmo,
\end{split}
\]
where
\[
    \frac{1}{r_2}=\frac{\abs{\mathfrak{J}_2}}{2}.
\]
Therefore, the same analysis as in \eqref{eq:FingMess} yields the result. 

If $\mathfrak{J}_2= \emptyset$, then $\mathfrak{J}_f\neq  \emptyset$. Therefore applying Fefferman--Stein's estimate \eqref{basic feffstein}, yields that \eqref{eq:peperoni} is bounded by 
\[
    \norm{f_1}_\bmo\norm{f_2}_\bmo\prod_{\mathfrak{J}_f}\norm{\Big(\sum_{j\geq 2k_0 } \,  \Big| (P_{j}^{u_\iota} \circ T^{\phase_\iota}_{d_\iota,j})(f_\iota)\Big|^2\Big)^{\frac{1}{2}}}_{L^{p_\iota}} \prod_{\mathfrak{J}_\infty} \norm{f_\iota}_\bmo 
\]
and using Lemma \ref{Lu-Be} concludes the discussion of this case.

\noindent{\bf The analysis concerning $d_0^\flat$.} 
The following lemma will be useful in to proving the desired regularity result.
\begin{lem}\label{handylem} Let $k_0$ be fixed, $0< p_0< \infty$ and $0<p_j\leq\infty$ so that 
\[
    \frac{1}{p_0}=\sum_{j=1}^N \frac{1}{p_j}.
\]
Then one has that
\[
    \sup_k \norm{P_{k_0}^0\brkt{G^U_k}}_{h^{p_0}}\lesssim c_{k_0}
    \prod_{j=1}^{N} \norm{f_j}_{X^{p_j}}.
\]
\end{lem}
\begin{proof} We shall give the proof for the case that $p_0\leq 1$. A small modification of the argument yields the case for $p_0>1$.

First we assume that $p_1+p_2<\infty$. We use the Littlewood-Paley characterisation of $h^{p_0}$, and the inclusion $\ell^{p_0} \subset \ell^1\subset \ell^2$. Then applying \cite{Triebel}*{p.17} and the fact that the frequency support of $\vartheta_j(D)P_{k_0}^0$ is included in a ball of radius $O(2^{k_0})$ followed by Hölder's inequality \eqref{eq:General_Holder} and Lemma \ref{lem:Wase-Lin}, we find that
\[
    \begin{split}
    \norm{P_{k_0}^0\brkt{G^U_k}}_{h^{p_0}} &\sim \norm{\brkt{\sum_{j=0}^{N(k_0)} \abs{\vartheta_j(D)P_{k_0}^0\brkt{G^U_k}}^2}^{1/2}}_{p_0}
    \\
    &\lesssim \brkt{ {{\sum_{j=0}^{N(k_0)} \norm{\vartheta_j(D) P_{k_0}^0\brkt{G^U_k}}_{p_0}^{{p_0}}}}}^{1/p_0}\lesssim \norm{{G^U_k}}_{p_0}\\
    &\leq 
    \prod_{\iota=1}^2{\norm{Q_{k}^{u_\iota} \circ T^{\phase_\iota}_{d_\iota}(f_\iota)}}_{h^{p_\iota}}
     \prod_{\iota\in \mathfrak{J}_2\cup \mathfrak{J}_f} \norm{P_{k}^{u_\iota} \circ T^{\phase_\iota}_{d_{\iota,k}}(f_\iota)}_{h^{p_\iota}}
     \prod_{\iota\in \mathfrak{J}_\infty} \norm{P_{k}^{u_\iota} \circ T^{\phase_\iota}_{d_{\iota,j}}(f_\iota)}_{L^\infty}\\
     &\lesssim \prod_{\iota=1}^2{\norm{f_\iota}}_{h^{p_\iota}}
     \prod_{\iota\in \mathfrak{J}_2\cup \mathfrak{J}_f} \norm{f_\iota}_{h^{p_\iota}}
     \prod_{\iota\in \mathfrak{J}_\infty} \norm{f_\iota}_{\bmo}
    \end{split}
\]

In the case that $p_1=\infty$ or $p_2=\infty$, a modification of the argument above, where one just uses that $\sup_{k\geq k_0} \norm{Q_k G}_{L^\infty} \lesssim \norm{G}_{\bmo},
$ yields the result. 
\end{proof}

Finally to deal with  $\sum_{k=k_0}^\infty  M_{\mathfrak m}\circ T^{\phase_0}_{d_0^\flat} \circ P_{k}^0\left[G^U_{k}\right](x)$ we observe that Lemma \ref{handylem} yields
\[
\begin{split}
         &\norm{\sum_{k=k_0}^\infty  2^{km_0}M_{\mathfrak m}\circ  T^{\phase_0}_{\Omega_0} \circ P_{k_0}^0\left[G^U_k\right]}^{p_0}_{h^{p_0}}\\
         &\leq {\sum_{k=k_0}^\infty  2^{km_0 {p_0}} \norm{P_{k_0}^0\left[G^U_k\right]}_{h^{p_0}}^{{p_0}}}\lesssim  \prod_{j=1}^{N} \norm{f_j}_{X^{p_j}}.
\end{split}
\]

Summing up and using the fact that $m_0<0$, we deduce the boundedness of $D(f_1, \dots, f_N)$, with target $h^{p_0}$.


\section{Space-time estimates for systems of dispersive PDEs}\label{space-time}
In this section we shall prove Theorem \ref{thm:PDE}, which amounts to 
showing Sobolev estimates for the solution $u$ of the system of coupled PDEs 
\begin{equation*}
\left\{ \begin{array}{l} i\partial_t u +  \varphi_0(D)\, u = T_{\zeta}\left( v_1,\dots, v_N\right)  \\
i\partial_t v_j +  \varphi_j(D)\, v_j = 0,\,\,\, j=1,\dots, N \\
\end{array} \right.
\quad \mbox{with} \quad
\left\{ \begin{array}{l} u(0,x) = 0  \\ v_j(0,x) = f_j (x), \,\,\, j=1,\dots, N, \end{array} \right. 
\end{equation*}
where $\varphi_j \in \mathcal{C}^{\infty}(\R^n \setminus 0)$, $f_j\in H^{\sigma_j, p_j}$, $\sigma_j\geq 0$, $j=0, \dots, \, N$ are assumed to be positively homogeneous of degree $s\in(0,\infty)$ and $T_\zeta$ is the multilinear multiplier given by \eqref{defn:tm} with symbol $\zeta\in S^{m_\zeta} (n, N)$ for some $m_ζ ≤ 0$, to be specified later. 
 The solution can be represented using the Duhamel formula as
\begin{equation}\label{representation of the solution}
u(t,x)= \int_{0}^{t}\int_{\R^{nN} } \zeta(\Xi)\, \prod_{j=1}^N \left(\widehat{f}_j(\xi_j)\, e^{ix\cdot\xi_j+ir \varphi_j(\xi_j)}\right)\,e^{i(t-r)\varphi_0( \xi_1+\cdots+\xi_N)}  \dd\Xi\, \dd r.
\end{equation}
This formula contains a multilinear oscillatory integral, and should therefore be suitable for analysis with the results of this paper. There are, however, two reasons why we cannot directly apply Theorems \ref{thm:main} and \ref{thm:main FIO}. Firstly, we must deal with the time dependency of $u$, and secondly, proving bounds in Sobolev spaces introduces more complicated amplitudes, which are a product of the multilinear amplitudes we have seen earlier and linear amplitudes in each variable. The following two results solve the first problem and extend regularity estimates of oscillatory integral operators with space-dependent phases to the corresponding time-dependent operators. We then proceed to prove Theorem \ref{thm:PDE} as a scholium to Theorems \ref{thm:main} and \ref{thm:main FIO}.

We shall start with the following lemma which yields time-dependent $L^p$ estimates for linear evolutions.
\begin{lem}\label{tidsuppskattning}
Let $\varphi\in \mathcal{C}^\infty(\R^n\setminus 0)$ be a phase function positively homogeneous of degree $s>0$. Then if $s\neq 1$ then for all $t ≥ 0$,
\begin{equation}\label{ litet t}
   \|\jap{ D}^{-sn|1/p - 1/2|}e^{it\varphi(D)}u\|_{L^p} \lesssim \jap{t}^{n|1/p - 1/2|} \|u\|_{L^p}, 
\end{equation}
and for $s=1$
\begin{equation}\label{ litet t}
   \|\jap{ D}^{-(n-1)|1/p - 1/2|}e^{it\varphi(D)}u\|_{L^p} \lesssim \jap{t}^{(n-1)|1/p - 1/2|} \| u\|_{L^p}. 
\end{equation}
\end{lem}
\begin{proof}
We only prove the case of $s\neq 1$, the remaining case is proven in a similar manner using Theorem \ref{linearhpthmfio}. 
First note that Theorem \ref{linearhpthmoio} yields that for $1<p<\infty$
\[
    \|e^{i\varphi(D)} \jap{ D}^{ -n s |\frac{1}{p}-\frac{1}{2}|}u\|_{L^p} \lesssim  \| u\|_{L^p}.
\]
To include $t$-dependence, we first note that in the case $t ≤ 1$,  $t\varphi(\xi)$ is a phase of order $s$ uniformly in $t$ and therefore satisfies the estimate
\begin{equation}\label{ litet t}
   \|e^{it\varphi(D)}u\|_{L^p} \lesssim \|\jap{ D}^{sn|1/p - 1/2|} u\|_{L^p} ≤ \jap{t}^a\|\jap{ D}^{sn|1/p - 1/2|} u\|_{L^p}, 
\end{equation}
for any $a ≥ 0$. When $t > 1$, we write $ m(p,s) =- ns|1/p - 1/2|$  and perform a change of variables (and using homogeneity of $\varphi$), finding
\begin{equation}\label{Akselsnille}
    \begin{split}
   ∫e^{ix·ξ + it\varphi(\xi)}\jap{ξ}^{m(p,s)}\hat u(ξ)\ddd ξ = t^{-n/s}∫e^{it^{-1/s}x·ξ + i\varphi(\xi)}\jap{t^{-1/s}ξ}^{m(p,s)}\hat u(t^{-1/s}ξ)\ddd ξ \\=
  t^{-m(p,s)/s}∫e^{it^{-1/s}x·ξ + i\varphi(\xi)}  σ_t(ξ)\widehat{ u(t^{1/s}\cdot)}(\xi)\ddd ξ,
\end{split}
\end{equation}
where $σ_t(ξ) = t^{m(p,s)/s}\jap{t^{-1/s}ξ}^{m(p,s)}$ satisfies $|∂_ξ^α σ_t(ξ)| ≤ C_α \jap{ξ}^{m(p,s) - |α|}$ when $t ≥ 1$ and $m(p,s) ≥ 0$. Therefore the $L^p-$bound given by Theorem \ref{linearhpthmoio}, \eqref{ litet t} and \eqref{Akselsnille} yield the desired result.
\end{proof}
A useful multilinear generalisation of this result is the following.
\begin{lem}\label{tidsuppskattningsmulti}
Let $\varphi_j\in \mathcal{C}^{\infty}(\R^{n} \setminus 0)$, $j=0,\ldots, N$, be phase functions that are homogeneous of degree $s>0$ and $\sigma\in S^{m}(n,N)$ with $m ∈ ℝ$.  Define
\[
    T^{(t)}_{\sigma} (f_1, \dots, f_N):=\int_{\R^{Nn}} e^{it\varphi_0( \xi_1+\cdots+\xi_N)} \sigma(\Xi)\, \prod_{j=1}^N \widehat{f}_j(\xi_j)\, e^{ix·\xi_j+it\varphi_j(\xi_j)}  \dd\Xi.
\]
Assume that for some $1<p_0, \ldots, p_N <\infty$ and $r_0, \ldots, r_N ∈ ℝ$ one has the estimate
\[
    \Vert \jap{D}^{-r_0}T^{(1)}_{\sigma} (f_1, \dots, f_N)\Vert_{L^p} ≤ C(\sigma, \Phi) \prod_{j=1}^{N}\Vert \jap D^{r_j}f_j\Vert_{L^{p_j}},
\]
 where $C(\sigma, \Phi)$ only depends on a finite number of seminorms of $\sigma$ and upper bounds on the size of a finite number of derivatives of $\varphi_j$. Then it follows that, for all $t\geq 0$
\[
    \Vert \jap D^{-r_0}T^{(t)}_{\sigma} (f_1, \dots, f_N)\Vert_{L^p} ≤ C(\sigma, \Phi) \langle{t}\rangle^{(\max(-m,0) + ∑_{j=0}^N \max(r_j,0))/s} \prod_{j=1}^{N}\Vert \jap D^{r_j}f_j\Vert_{L^{p_j}}.
\]
\end{lem}
\begin{proof}
For $0 ≤ t\leq 1$, there is an upper bound on the derivatives of $t\Phi$ that is uniform in $t$, so this case is clear. When $t > 1$, we let $g_j = \jap D^{r_j}f_j$, so that 
\begin{align*}
    \jap D^{-r_0}T^{(t)}_{σ}(&f_1, \ldots, f_N)(x) =\\[.4em]
    &\phantom{={}} ∫_{ℝ^{Nn}} e^{ix·(ξ_1 +\cdots +ξ_N) + itφ_0(ξ_1 + \cdots + ξ_N) + ∑_{j=1}^N itφ_j(ξ_j)}\\
    &\qquad×\ \jap{ξ_1 + \cdots + ξ_N}^{-r_0}σ(Ξ)∏_{j=1}^N \jap{ξ_j}^{-r_j}\hat g_j(ξ_j)\dd Ξ\\
    &= t^{-Nn/s} ∫_{ℝ^{Nn}} e^{it^{-1/s}x·(ξ_1 +\cdots +ξ_N) + iφ_0(ξ_1 + \cdots + ξ_N) + ∑_{j=1}^N iφ_j(ξ_j)}\\
    &\qquad×\ \jap{t^{-1/s}(ξ_1 + \cdots + ξ_N)}^{-r_0}σ(t^{-1/s}Ξ)∏_{j=1}^N \jap{t^{-1/s}ξ_j}^{-r_j}\hat g_j(t^{-1/s}ξ_j)\dd Ξ\\
    &= t^{(\max(-m,0) + ∑_{j=0}^N \max(r_j,0) - Nn)/s} ∫_{ℝ^{Nn}} e^{it^{-1/s}x·(ξ_1 +\cdots +ξ_N) + iφ_0(ξ_1 + \cdots + ξ_N) + ∑_{j=1}^N iφ_j(ξ_j)}\\
    &\qquad×\ \jap{ξ_1+ \cdots + ξ_N}^{-r_0}t^{\min(m,0)/s}σ(t^{-1/s}Ξ)\bigg(∏_{j=1}^N \jap{ξ_j}^{-r_j}\hat g_j(t^{-1/s}ξ_j)\bigg)\\
    &\qquad×\ \frac{t^{-\max(r_0,0)/s}\jap{t^{-1/s}(ξ_1 + \cdots + ξ_N)}^{-r_0}}{\jap{ξ_1 + \cdots + ξ_N}^{-r_0}}∏_{j=1}^N \frac{t^{-\max(r_j,0)/s}\jap{t^{-1/s}ξ_j}^{-r_j}}{\jap{ξ_j}^{-r_j}}\dd Ξ\\
    &= t^{(\max(-m,0) + ∑_{j=0}^N \max(r_j,0))/s}\\
    &\quad×\ S_0 \jap D^{-r_0}T^{(1)}_{σ_t}(\jap D^{-r_1}S_1 g_1(t^{1/s}·),\,\ldots,\,\jap D^{-r_N}S_N g_N(t^{1/s}·))(t^{-1/s}x),
\end{align*}
where
\begin{align*}
    S_j &= t^{-\max(r_j,0)/s}\jap{t^{-1/s}D}^{-r_j}\jap{D}^{r_j},\\
    σ_t(Ξ) &= t^{\min(m,0)/s}σ(t^{-1/s}Ξ).
\end{align*}
Now, $σ_t ∈ S^m(n,N)$ uniformly in $t$, so we can use the known boundedness of $\jap{D}^{-r_0}T^{(1)}_{σ_t}$. The operators $S_j$ are furthermore Mikhlin multipliers uniformly in $t$ and hence bounded $L^p → L^p$. It follows that
\[
    \|\jap{D}^{-r_0}T^{(t)}_{σ}(f_1,\ldots,f_N)\|_{L^{p_0}} \lesssim t^{(\max(-m,0) + ∑_{j=0}^N \max(r_j,0))/s}∏_{j=1}^N\|\jap D^{r_j}f_j\|_{L^{p_j}}.\qedhere
\]
\end{proof}

Now let us return to the Duhamel representation \eqref{representation of the solution}.
Here we set $$T^{(r)}_{\zeta} (f_1, \dots, f_N)(x):=\int_{\R^{nN} } \zeta(\Xi)\, \prod_{j=1}^N \left(\widehat{f}_j(\xi_j)\, e^{ix\cdot\xi_j+ir \varphi_j(\xi_j)}\right)  \dd\Xi,$$ 
and observe that
\[
    u(t,x)= \int_{0}^{t} e^{i(t-r)\varphi_0(D)}  \jap{D}^{m(p_0,s)} \jap{D}^{-m(p_0,s)}T^{(r)}_{\zeta} (f_1, \dots, f_N)(x) \dd r.
\]
Let 
\[
    \sigma_0=\varkappa+m_c-m_\zeta,\qquad \varkappa:= \min_{j=1,\ldots,N} \sigma_j,
\]
where $m_c = m_c(s)$ is as in the statement of Theorem \ref{thm:PDE}.
From this and Lemma \ref{tidsuppskattning} we immediately obtain
\begin{equation}\label{tidsintegralen}
    \Vert u\Vert_{H^{\sigma_0,p_0}}\lesssim \int_0^t  \jap{t-r}^{-m(p_0, s)/s}\Vert \jap{D}^{-m(p_0,s)}T^{(r)}_{\zeta} (f_1, \dots, f_N)\Vert_{H^{\sigma_0,p_0}} \dd r,
\end{equation}
for $1<p_0<\infty$.
Using Lemma \ref{tidsuppskattningsmulti} it will therefore be enough for us to study the right-hand norm in the case where $r=1$. Now, using the decomposition of Section \ref{sec:freqdecom} we can decompose $\zeta$ and reduce the analysis of $T^{(1)}_{\zeta}$ to the study of multilinear operators $T_{\zeta_{0}}$, $T_{\zeta_{1}}$ and $T_{\zeta_{1,2}}$. It should however be noted that for these terms the method only takes advantage of the added regularity on the first argument  (i.\,e.\@{} $f_1$). For the similar terms $T_{\zeta_2}$, $T_{\zeta_{2,3}}$, etc.\@{} one can take advantage of a different $σ_j$. This is the reason why $\varkappa$ is the minimum of the $σ_j$, $j = 1, \ldots, N$.

\subsubsection*{Treatment of $T_{ζ_0}$} 
Here we make use of the representation given in (\ref{final form of sigma0}), which in our case with $x$-independent amplitude translates to
\begin{align*}
    &\textbf I := \jap{D}^{-m(p_0,s)}T_{ζ_0} (f_1,\dots,f_N) \\[.2em]
    &= \sum_{K\in\Z^{nN}} a_K\jap{D}^{-m(p_0,s)}\theta(D/\sqrt{N})\bigg(\prod_{j=1}^N T_{\theta}^{\phase_j}\circ\tau_{\frac{2\pi k_j}{L}}(f_j)\bigg).
\end{align*}
The method in Subsection \ref{boundedness of Tsigma0} can then be carried out to show that for any $\sigma\in \R$
\[
    \Vert {\bf{I}}\Vert_{H^{\sigma, p_0}}\lesssim  \prod_{j=1}^{N} \Vert f_j\Vert_{X^{p_j}} \lesssim ∏_{j=1}^N \|f_j\|_{H^{σ_j,p_j}}.
\]

\subsubsection*{Treatment of $T_{ζ_1}$}

Using \eqref{t_piecedecomp for sigma1} and \eqref{the magic operator}, with the same notation as was introduced there, its $L^p$-boundedness can be inferred from that of the multilinear operator
\begin{align*}\label{basicsystem estim1}
     &\textbf{II}:= \jap D^{-m(p_0,s)}T_{ζ_1}(f_1,\ldots,f_N) = ∫\widetilde{\textbf{II}}_U\, \frac1{(1+|U|^2)^M}  \, \ddd U,
\end{align*}     
where
\begin{align*}
    \widetilde{\textbf{II}}_U &=\jap{D}^{-m(p_0,s)} { \sum_{k\geq k_0 } \chi_0(2D)\,Q_{k}^0  \left[(Q_{k}^{u_1} \circ T^{ \phase_1}_{b_1})(f_1)\prod_{j=2}^N(P_{k}^{u_j} \circ T^{\phase_j}_{b_{j,k}})(f_j)\right]}\\
    &=        { \sum_{k\geq k_0 } \chi_0(2D)    \jap{D}^{-m(p_0,s)}\circ \,Q_{k}^0  \left[(Q_{k}^{u_1} \circ T^{\phase_1}_{b_1 |\cdot|^{-\sigma_1}})(|D|^{\sigma_1}f_1)\prod_{j=2}^N(P_{k}^{u_j} \circ T^{\phase_j}_{b_{j,k} })(f_j)\right]}\\
        &=        \chi_0(2D)    \jap{D}^{-m(p_0,s)}\abs{D}^{m_\zeta+m(p_0,s)-m_c(s)-\sigma_1}\\ 
        &\ \circ \sum_{k\geq k_0 }  \,Q_{k}^1  \left[(Q_{k}^2 \circ T^{\phase_1}_{\widetilde{b_1} })(\abs{D}^{\sigma_1}\chi_0(2D)f_1)\prod_{j=2}^N(P_{k}^{u_j} \circ T^{\phase_j}_{\widetilde{b_{j,k}} })(f_j)\right],
 \end{align*}
where $b_1\in S^{m_1}$ and $b_{j,k}\in S^{m_j}$. 
$Q^1_k$ has symbol
\[
    {\phi}_k(\xi) |2^{-k}\xi|^{-m_\zeta+m_c(s){-m_0(p_0,s)}+\sigma_1} ,
\]
$Q^2_k$ has symbol
\[
    {\psi}_k(\xi) |2^{-k}\xi|^{m_ζ-m(p_1,s)-\sigma_1} e^{2^{-k}\xi\cdot u_1},
\]
and we define 
\begin{align*}
    \widetilde{b_1}(\xi)&=\chi_0(\xi)\abs{\xi}^{m(p_1,s)}\in S^{m(p_1,s)},\\
    \widetilde{b_{j,k}}(\xi)&=2^{m(p_j,s)k}\omega_k(\xi), \qquad j=2,\ldots, N.
\end{align*}
Now since the operator $\sum_{k\geq k_0 }  \,Q_{k}^1  \left[(Q_{k}^2 \circ T^{\phase_1}_{\widetilde{b_1} })( f_1)\prod_{j=2}^N(P_{k}^{u_j} \circ T^{\phase_j}_{b_{j,k} })(f_j)\right]$ is of the form (\ref{defn:thick I}) the boundedness of the latter yields that
\[
    \Vert {\bf{II}}\Vert_{H^{\sigma_0,p_0}}\lesssim \Vert {\bf{II}}\Vert_{H^{\sigma_1+m_c-m_\zeta,p_0}}\lesssim \Vert f_1\Vert_{H^{\sigma_1, p_1}} \prod_{j=2}^{N} \Vert f_j\Vert_{X^{p_j}}
    \lesssim ∏_{j=1}^N\|f_j\|_{H^{σ_j,p_j}}.
\]

\subsubsection*{Treatment of $T_{ζ_{1,2}}$} For this part we will need to invoke an interpolation argument. 
To that end, we fix $\sigma_j$'s , $0\leq j\leq N$ with $\sigma_0 ≥ 0$. Then the goal is to show that the $N$-linear operator $W$ given by
\[
    W(f_1, \dots, f_N):= \jap{D}^{\sigma_0-m(p_{0}, s)}T_{\zeta_{1,2}} (\jap{D}^{-σ_1} f_1, \dots, \jap D^{-σ_N}f_N)
\]
is bounded $\prod_{j}X^{p_j} \to X^{p_0}$, provided that $m_\zeta = m_c - \sigma_0 + \varkappa$. Observe that $m_c$ depends linearly on the $1/p_j$ between any two adjacent endpoints in Lemma \ref{cor:endpointcases}, and hence the same goes for $m_ζ$.
We can therefore use the interpolation argument in Lemma \ref{cor:endpointcases} on $W$.

Just as in the treatment of \textbf{II}, we only use the Sobolev regularity in $f_1$, and that of $f_2,\ldots, f_N$ will only be used in the analogous estimates for other $\zeta_{i,j}$.


Now, using the representation \eqref{t_piecedecomp for sigma N+1 bilin}, we need to study the boundedness of 
\begin{align*}\label{basicsystem estim2}
    &\textbf{III} := \jap D^{-m(p_0,s)}T_{ζ_{1,2}}(f_1,\ldots,f_N) = ∫\widetilde{\textbf{III}}_U\, \frac1{(1+|U|^2)^M}  \, \ddd U,
\end{align*}
where $U = (u_1,\dots,u_N)$ and
\begin{align*}
    \widetilde{\textbf{III}}_U &=\sum_{k\geq k_0}^\infty  M_{\mathfrak m} d_0(D)\jap{D}^{-m(p_0,s)}  P_{k}^0\left[(Q_{k}^{u_1} \circ T^{\phase_1}_{d_1})(f_1)\, (Q_{k}^{u_2} \circ T^{\phase_2}_{d_2})(f_2)\,\prod_{j=3}^N(P_{k}^{u_j} \circ T^{\phase_j}_{d_{j,k}})(f_j)\right],
\end{align*}
with $M_{\mathfrak m}$ being the operation of multiplication by $\mathfrak m(k, U)$, which is uniformly bounded in $k$. Moreover $d_0=2^{k(m_\zeta-m_c+m(p_0,s))} \omega_k(\xi)$ and the amplitudes for each OIO are defined by \eqref{the d's}.

We shall consider the norm of $\widetilde{\textbf{III}}_U$ in $H^{\sigma_0,p_0}$ where $p_0=2$, $p_0=1$ and $p_0=\infty$, which by duality corresponds to estimating 
\[
     S := ∫\widetilde{\textbf{III}}_U(x)f_0(x)\,\dd x
\]
with $f_0\in H^{-\sigma_0,p'_{0}}$ ($p'_{0}$ is the H\"older dual of $p_0$).
First we observe that 
\[
    P^{0}_{k}= P_{k_{0}} + \sum_{\ell=k_{0}+1}^{k} Q_{\ell}
\]
and therefore one can write $S = S_\textsc p + S_\textsc q$ with
\begin{align*}
    S_\textsc p &= ∫\sum_{k\geq k_0}^\infty M_{\mathfrak m} d_0(D)\jap{D}^{-m(p_0,s)} P_{k_0}f_0(x)\, Q_{k}^{u_1} T^{\phase_1}_{d_1}f_1(x)\, Q_{k}^{u_2}T^{\phase_2}_{d_2}f_2(x)\,\prod_{j=3}^N P_{k}^{u_j} T^{\phase_j}_{d_{j,k}}f_j(x)\,\dd x\\
    S_\textsc q &= ∫\sum_{k\geq k_0}^\infty\sum_{\ell=k_{0}+1}^{k} M_{\mathfrak m} Q_ℓ d_0(D)\jap{D}^{-m(p_0,s)} f_0(x)\, Q_{k}^{u_1} T^{\phase_1}_{d_1}f_1(x)\, Q_{k}^{u_2}T^{\phase_2}_{d_2}f_2(x)\,\prod_{j=3}^N P_{k}^{u_j} T^{\phase_j}_{d_{j,k}}f_j(x)\,\dd x
\end{align*}

To show the needed boundedness of these parts, we shall rely on the method laid out in detail in Section 8.1 of \cite{RRS2}. The terms $S$, and $S_\textsc p$ correspond in that text to the expressions (60) and (61), respectively. For the term $S_\textsc q$ we 
note that using the condition 
$\varkappa + m_c - m_{\zeta}\geq 0$, we have that
\begin{multline}\label{L2forumsigma}
\sum_{k \geqslant k_{0}}\sum_{\ell=k_0+1}^{k} \Big\vert \int \Big(  M_{\mathfrak m} Q_{\ell}  d_0(D)\jap{D}^{-m(p_0,s)} f_{0}\Big)(x)\left(Q_{k}^{u_{1}} \circ T_{d_{1}}^{\varphi_{1}}\right)\left(f_{1}\right)(x) \\ \times\left(Q_{k}^{u_{2}} \circ T_{d_{2}}^{\varphi_{2}}\right)\left(f_{2}\right)(x) \prod_{j=3}^{N}\left(P_{k}^{u_{j}} \circ T_{d_{j,k}}^{\varphi_{j}}\right)\left(f_{j}\right)(x) \mathrm{d} x\Big \vert\\ \leq
\sum_{k \geqslant k_{0}}\sum_{\ell=k_0+1}^{k} 2^{(k-\ell)(m(p_0,s)-\sigma_1 -m_c +m_\zeta)} \Big\vert \int \Big(  M_{\mathfrak m}2^{-km(p_0,s)} Q_{\ell}  \tilde{d}_0(D) \jap{D}^{-\sigma_1 -m_c +m_\zeta}f_{0}\Big)(x)\\ \times\left(Q_{k}^{u_{1}} \circ T_{d_{1}}^{\varphi_{1}}\right)\left(|D|^{\sigma_1}f_{1}\right)(x) \left(Q_{k}^{u_{2}} \circ T_{d_{2}}^{\varphi_{2}}\right)\left(f_{2}\right)(x) \prod_{j=3}^{N}\left(P_{k}^{u_{j}} \circ T_{d_{j,k}}^{\varphi_{j}}\right)\left(f_{j}\right)(x)\, \mathrm{d} x\Big \vert
\\
\leq
\sum_{k \geqslant k_{0}}\sum_{\ell=k_0+1}^{k} \Big\vert \int \Big(  M_{\mathfrak m} Q_{\ell}  \tilde{d}_0(D) \jap{D}^{-\sigma_1 -m_c +m_\zeta}f_{0}\Big)(x)\left(Q_{k}^{u_{1}} \circ T_{d_{1}}^{\varphi_{1}}\right)\left(|D|^{\sigma_1}f_{1}\right)(x) \\ \times\left(Q_{k}^{u_{2}} \circ T_{d_{2}}^{\varphi_{2}}\right)\left(f_{2}\right)(x) \prod_{j=3}^{N}\left(P_{k}^{u_{j}} \circ T_{d_{j,k}}^{\varphi_{j}}\right)\left(f_{j}\right)(x) \mathrm{d} x\Big \vert
 \\
=
\sum_{\ell=k_0}^{\infty}\sum_{k=0}^{\infty} \Big\vert \int \Big(  M_{\mathfrak m} Q_{\ell}  \tilde{d}_0(D) \jap{D}^{\varkappa-\sigma_0-\sigma_1}f_{0}\Big)(x)\left(Q_{k+\ell}^{u_{1}} \circ T_{d_{1}}^{\varphi_{1}}\right)\left(|D|^{\sigma_1}f_{1}\right)(x) \\
\times\left(Q_{k+\ell}^{u_{2}} \circ T_{d_{2}}^{\varphi_{2}}\right)\left(f_{2}\right)(x) \prod_{j=3}^{N}\left(P_{k+\ell}^{u_{j}} \circ T_{d_{j,k}}^{\varphi_{j}}\right)\left(f_{j}\right)(x)\, \mathrm{d} x\Big \vert
\end{multline}
where $\tilde{d}_0(\xi)=  \omega_k(\xi)$. This last expression corresponds to the sum in $k$ of expression (62) in \cite{RRS2}.

With this set, one can follow the procedure in \cite{RRS2} to show the required end-point estimates. However, in order for every step of that proof to translate to this setting, we need to show some additional facts about our terms. 

First we consider the target space $H^{\sigma_0, 2},$ and to make use of duality take $f_0$ such that $\jap{D}^{-\sigma_0}f_0\in L^2$. 
To deal with 
$S_\textsc p$ we hence have to estimate
\begin{multline*}\label{eq:lowfreq11}
\sum_{k\geq k_0}\Big|
\int \Big(M_{\mathfrak m} P_{k_0} d_0(D)\jap{D}^{-m(p_0,s)} f_0\Big)(x)\\ \times(Q_{k}^{u_1} \circ T^{\phase_1}_{d_1})(f_1)(x)\, (Q_{k}^{u_2} \circ T^{\phase_2}_{d_2})(f_2)(x)\,\prod_{j=3}^N(P_{k}^{u_j} \circ T^{\phase_j}_{d_j})(f_j)(x)\,\dd x\Big|
\end{multline*}
Now since $k_0$ is fixed, the symbol of the multiplier $P_{k_0}$ is a Schwartz function and therefore 
\begin{equation*}\label{kernel-maximal}
\begin{aligned}
    &\quad M_{\mathfrak m} P_{k_0} d_0(D)\jap{D}^{-m(p_0,s)} f_0  \\
    &= (P_{k_0} \jap{D}^{-m(p_0,s)} \jap{D}^{\sigma_1+m_c -m_\zeta}\circ d_0 (D) \circ P_k \circ M_{\mathfrak m}) (\jap{D}^{-\sigma_1-m_c +m_\zeta}f_0)\\
    &= K * ((P_k \circ M_{\mathfrak m}) (\jap{D}^{-\sigma_1-m_c +m_\zeta}f_0)),  
\end{aligned}
\end{equation*}
for $k \geq k_0$, with $|K(\cdot)| \lesssim \langle \cdot \rangle^{-N},$ for any $N\geq 0$, which shows that this term has the required form for the steps on page~36 in \cite{RRS2} to go through.

Following those steps, we therefore see that \textbf{III} is bounded in $L^2$ provided that for $j=1,2$ and $f\in \bmo$ the measure
\[
    \mathrm{d} \mu_{k}(x, t)=\sum_{\ell=0}^{\infty}\left|\left(Q_{k+\ell}^{u_{j}} \circ T_{d_{j}}^{\varphi_{j}}\right)(f)(x)\right|^{2} \delta_{2^{-\ell}}(t) \mathrm{d} x
\]
is Carleson with a decay in $ℓ$ in the Carleson norm. However, in Proposition~\ref{lem:smallcarlnorm} it was shown that the Carleson norm is bounded by a multiple of $2^{-εk}\|f\|_\bmo^2$, for some $ε > 0$.  This decay in $k$ is needed to be able to deal with double sum in \eqref{L2forumsigma}  in various cases that are handled below.

This fact enables us to use the arguments in Section 8.1 on page 35 of \cite{RRS2}, in accordance with case (\textit{ii}) of Lemma \ref{cor:endpointcases} to conclude that 
\[
    \Vert \jap{D}^{\sigma_0}{\bf{III}}\Vert_{L^2}\lesssim \Vert\jap{D}^{\sigma_1}f_1\Vert_{X^{p_1}}   \prod_{j=2}^{N} \Vert f_j\Vert_{X^{p_j}} \lesssim ∏_{j=1}^N\|\jap D^{σ_j}f_j\|_{X^{p_j}},
\]
and therefore
\[
    \|W(f_1,\ldots,f_N)\|_{L^2} \lesssim ∏_{j=1}^N \|f_j\|_{X^{p_j}}.
\]

Next, we deal with the target space with norm $\|\jap D^{σ_0}· \|_{h^1}$, and therefore take $f_0$ such that $\jap{D}^{-\sigma_0}f_0\in \bmo$.  Therefore if $\jap{D}^{\sigma_1}f_1\in \bmo$ and $f_j\in\bmo$ for $j=2, \dots, M$, then for any $3\leq M\leq N$ it is not hard (mainly using Proposition \ref{lem:smallcarlnorm}) to see that the measure
\[
\begin{aligned} \mathrm{d} \mu_k(x, t) &:=\sum_{\ell=0}^{\infty} \Big(Q_{\ell}M_{\mathfrak m}   \tilde{d}_0(D) \jap{D}^{-\sigma_1 -m_c +m_\zeta}f_{0}(x)\Big) \\ & \times\left[\left(Q_{k+\ell}^{u_{1}} \circ T_{d_{1}}^{\varphi_{1}}\right)\left(|D|^{\sigma_1}f_{1}\right)(x)\left(Q_{k+\ell}^{u_{1}} \circ T_{d_{2}}^{\varphi_{2}}\right)\left(f_{2}\right)(x) \prod_{j=3}^{M}\left(P_{k+\ell}^{u_{j}} \circ T_{d_{j}}^{\varphi_{j}}\right)\left(f_{j}\right)(x)\right] \mathrm{d} x\, \delta_{2^{-l}}(t) \end{aligned}
\]
is a Carleson measure with the Carleson norm bounded by a multiple of
\[
    2^{-εk} \norm{\jap{D}^{-\sigma_0}f_0}_\bmo \norm{\jap{D}^{\sigma_1}f_1}_\bmo \prod_{j=2}^{M} \Vert f_j\Vert_{\bmo},
\]
 for some $ε > 0$.
Moreover by estimate \eqref{eq:uniform_bounddness} we also have that
\[
  \sup _{\ell \geqslant k_{0}}\left\|Q_{\ell}M_{\mathfrak m}   \tilde{d}_0(D) \jap{D}^{-\sigma_1 -m_c +m_\zeta}f_{0}\right\|_{L^{\infty}} \lesssim\left\|\jap{D}^{-\sigma_0}f_{0}\right\|_{\text {bmo }} 
\]
and 
\[
    \sup _{\ell \geqslant k_{0}}\left\|\left(Q_{k+\ell}^{u_{j}} \circ T_{d_{j}}^{\varphi_{j}}\right)\left(f_{j}\right)\right\|_{L^{\infty}} \lesssim\left\|f_{j}\right\|_{\mathrm{bmo}} \quad \text { for } j=1,2 \text { when } p_{j}=\infty,
\]
where the hidden constant in the above estimate is uniform in $k$.
These facts together with estimates \eqref{ineq:carll2}, \eqref{ineq:carlh1} and \eqref{ineq:quadraticestimate} enable us to run the arguments of Section 8.2 on page 40 of \cite{RRS2} to prove various boundedness results corresponding to the cases {$($\it{iii}$)$} and {$($\it{iv}$)$} of Lemma \ref{cor:endpointcases} and finally arrive at
\[
    \Vert \jap{D}^{\sigma_0} {\bf{III}}\Vert_{h^1}\lesssim \Vert \jap{D}^{\sigma_1}f_1\Vert_{X^{p_1}}  \prod_{j=2}^{N} \Vert f_j\Vert_{X^{p_j}} \lesssim ∏_{j=1}^N \|\jap D^{σ_j}f_j\|_{X^{p_j}},
\]
and hence
\[
    \|W(f_1,\ldots,f_N)\|_{h^1} \lesssim ∏_{j=1}^N \|f_j\|_{X^{p_j}}.
\]

The last case to deal with is when $f_0$ in the duality arguments above has the property that $\jap{D}^{-\sigma_0}f_0\in h^1$.
Here we observe that the measure
\[
    \mathrm{d} \mu_k(x, t) =\sum_{\ell=0}^{\infty}\left(Q_{k+\ell}^{u_{1}} \circ T_{d_{1}}^{\varphi_{1}}\right)\left(|D|^{\sigma_1}f_{1}\right)(x)\left(Q_{k+\ell}^{u_{1}} \circ T_{d_{2}}^{\varphi_{2}}\right)\left(f_{2}\right)(x) \prod_{j=3}^{M}\left(P_{k+\ell}^{u_{j}} \circ T_{d_{j}}^{\varphi_{j}}\right)\left(f_{j}\right)(x)\, \mathrm{d} x \delta_{2^{-l}}(t) 
\]
is Carleson with Carleson norm bounded by a multiple of
\[
   2^{-εk}  \norm{\jap{D}^{\sigma_1}f_1}_\bmo \prod_{j=2}^{N} \Vert f_j\Vert_{\bmo},
\]
 for some $ε > 0$. Therefore \eqref{ineq:carlh1} yields that 
\[
    \Vert \jap{D}^{\sigma_0}{\bf{III}}\Vert_{\bmo}\lesssim \Vert\jap{D}^{\sigma_1}f_1\Vert_{X^{p_1}}   \prod_{j=2}^{N} \Vert f_j\Vert_{X^{p_j}} \lesssim ∏_{j=1}^N\|\jap D^{σ_j}f_j\|_{X^{p_j}},
\]
yielding
\[
    \|W(f_1,\ldots,f_N)\|_\bmo \lesssim ∏_{j=1}^N \|f_j\|_{X^{p_j}}.
\]

With all the end point estimates set, we can by interpolation finally deduce that 
\[
    \|W(f_1,\ldots,f_N)\|_{X^{p_0}} \lesssim ∏_{j=1}^N\|f_j\|_{X^{p_j}},\quad p_0,\ldots,p_N ∈ [0,∞],
\]
which means that
\[
    \|\textbf{III}\|_{H^{σ_0,p_0}} \lesssim  ∏_{j=1}^N \|f_j\|_{H^{σ_j,p_j}}, \quad p_0,\ldots, p_N ∈ (0,∞).
\]

Returning now to \eqref{tidsintegralen}, we recall that $T^{(1)}_{\zeta}$ is a sum of operators, of the type $T_{\zeta_0}$, $T_{\zeta_1}$ and $T_{\zeta_{1,2}}$ and the bounds obtained above for \textbf{I}, \textbf{II} and \textbf{III} can therefore be used to show that 
\begin{equation*}
    \Vert \jap{D}^{-m(p_0,s)} T^{(1)}_{\zeta}(f_1, \dots, f_N)\Vert_{H^{\sigma_0, p_0}}\lesssim \prod_{j=1}^{N}\Vert  f_j\Vert_{H^{\sigma_j, p_j}}.
\end{equation*}
Lemma \ref{tidsuppskattningsmulti} then yields that
\[
    \Vert \jap{D}^{-m(p_0,s)+\varkappa+m_c-m_\zeta} T^{(r)}_{\zeta}(f_1, \dots, f_N)\Vert_{L^{ p_0}}\lesssim \langle{r}\rangle^{(-m_\zeta +\varkappa'+ ∑_{j=1}^N \sigma_j)/s} \prod_{j=1}^{N}\Vert  f_j\Vert_{H^{\sigma_j, p_j}},
\]
where $\varkappa':= \max(m(p_0,s)-\varkappa-m_c+m_\zeta, 0).$
Thus we conclude that for the solution $u$ in \eqref{tidsintegralen} one has
\begin{equation*}\label{final PDE 1}
    \Vert u(t, \cdot)\Vert_{H^{\sigma_0,p_0}(\R ^n)}\lesssim \int_{0}^{t} \jap{t-r}^{-m(p_0, s)/s} \langle{r}\rangle^{(-m_\zeta +\varkappa'+ ∑_{j=1}^N \sigma_j)/s} \dd r\, \prod_{j=1}^{N}\Vert  f_j\Vert_{H^{\sigma_j, p_j}}
\end{equation*}
from which one obtains the space-time estimate
\[
   \Vert u \Vert_{L^{q}([0,T])\,H^{\varkappa+m_c-m_{\zeta},p_0}(\R ^n)}\leq C_{T} \prod_{j=1}^{N}\Vert  f_j\Vert_{H^{\sigma_j, p_j}},
\]
which is valid for any $q\in[1,\infty]$, any $T\in (0,\infty).$ Theorem \ref{thm:PDE} is thereby proven.

\end{document}